% On Subgroups Of Coxeter Groups by Warren Dicks and Ian Leary
% A plain TeX file

%\input BoxedEPS.tex
%\SetTexturesEPSFSpecial
%\HideDisplacementBoxes

\input psfig.sty

\font\Bbten=msbm10
\font\Bbseven=msbm7

\magnification=\magstep 1
\pageno=1

\def\BBb#1{\hbox{\Bbten #1}}
\def\Bbb#1{\hbox{\Bbseven #1}}

\def\Zz{{\BBb Z}}
\def\zz{{\Bbb Z}}
\def\Qq{{\BBb Q}}

\def\Ff{{\BBb F}}
\def\ff{{\Bbb F}}
\def\Cc{{\BBb C}}

% The next line gets rid of things like  CW-
% complex.
\exhyphenpenalty=10000

% This is the same as the usual beginsection, except that it doesn't
% mind leaving up to 10% of a page blank, whereas the original doesn't
% mind leaving up to 30% of a page blank.
\def\beginsection#1\par{\vskip0pt plus.1\vsize\penalty-250
\vskip0pt plus-.1\vsize\bigskip\vskip\parskip
\message{#1}
\leftline{\bf#1}\nobreak\smallskip\noindent}

\def\printpictures{
\pageinsert
%\null \vskip 1truecm
%\centerline{\BoxedEPSF{Lyndon scaled 850}}
%\centerline{
\vskip -2in
\psfig{figure=lyndon,height=6in,bbllx=-100pt,bblly=-150pt,bburx=400pt,bbury=350pt}
%} 
\endinsert
\pageinsert
%\null \vskip 4truecm
%\centerline{\BoxedEPSF{Obverse scaled 850}} 
%\centerline{
\psfig{figure=obverse,height=6in,bbllx=-100pt,bblly=-300pt,bburx=400pt,bbury=200pt}
%}
\endinsert}

\countdef\picturesdone=42

\picturesdone=0

\def\picturesout{
\goodbreak     %needed in case picturesout is called at the bottom of
%               an odd page.
\ifnum\picturesdone=0
      \immediate\ifodd\pageno
           \printpictures
           \message{Pictures appear on pages \folio+1 and \folio+2.}
           \picturesdone=1
      \fi
\fi}

\def\sqr#1#2{{\vcenter{\vbox{\hrule height.#2pt
   \hbox{\vrule width.#2pt height#1pt \kern#1pt
     \vrule width.#2pt}
   \hrule height.#2pt}}}}
\def\whis{\mathchoice\sqr34\sqr34\sqr{2.1}3\sqr{1.5}3}
\def\qed{\ifmmode \eqno{\whis}\else
{{\unskip\nobreak\hfil\penalty50\hskip 2em\hbox{}\nobreak\hfil
${\whis}$\parfillskip=0pt\par\medskip}}\fi}

\def\inv #1{{\bar #1}}
\def\sinv{\inv s}
\def\tinv{\inv t}
\def\uinv{\inv u}
\def\vinv{\inv v}
\def\winv{\inv w}
\def\xinv{\inv x}
\def\yinv{\inv y}
\def\zinv{\inv z}

\def\cd{{\rm cd}}
\def\vcd{{\rm vcd}}
\def\Hom{{\rm Hom}}
\def\Ext{{\rm Ext}}
\def\Tor{{\rm Tor}}
\def\Aut{{\rm Aut}}
\def\sint{{\rm int}}
\def\containedin{\subseteq}
\def\empty{\emptyset}
\def\remark{\par\noindent{\bf Remark.} }
\def\definition{\par\noindent{\bf Definition.} }
\def\pf{\par\noindent{\bf Proof.} }
\def\proof{\pf}
\def\pra{\par}
\def\mapright#1{\smash{\mathop{\longrightarrow}\limits^{#1}}}
\def\mapdown#1{\Big\downarrow\rlap{$\vcenter{\hbox{$\scriptstyle#1$}}$}}

\def\bes{1}
\def\bie{2}
\def\bre{3}
\def\brp{4}
\def\bro{5}
\def \chs{6}
\def\chis{7}
\def\dav{8}
\def\daj{9}
\def\did{10}
\def\dil{11}
\def\gre{12}
\def\ham{13}
\def\hst{14}
\def\mun{15}
\def\rob{16}
\def\gap{17}

\centerline{\bf On subgroups of Coxeter groups}
\centerline{by}
$$\hbox{\vbox{\hbox{Warren Dicks}\hbox{Departament de Matem\`atiques,}
\hbox{Universitat Aut\`onoma de Barcelona,}
\hbox{E 08193, Bellaterra (Barcelona),}\hbox{Spain}
}\qquad\vbox{\hbox{I. J. Leary}\hbox{Faculty of Mathematical Studies,}
\hbox{University of Southampton,}\hbox{Southampton SO17 1BJ,}
\hbox{England}}} $$
\centerline {\it Dedicated to the memory of Brian Hartley}

\beginsection 1. Outline.

For a finitely generated
Coxeter group $\Gamma$,  its virtual cohomological dimension
over a (non-zero, associative) ring $R$, denoted $\vcd_R\Gamma$,  is
finite and has been described [\dav,\bes,\dil,\ham].  In [\dav], M.~Davis
introduced a contractible $\Gamma$-simplicial complex with finite
stabilisers.  The dimension of such a complex gives an upper bound for
$\vcd_R\Gamma$.  In [\bes], M.~Bestvina gave an algorithm for
constructing an $R$-acyclic $\Gamma$-simplicial complex with finite
stabilisers of dimension exactly $\vcd_R\Gamma$, for $R$ the integers or
a prime field;  he used this to exhibit a group whose  cohomological
dimension over the integers  is finite but strictly greater than its
cohomological dimension  over the rationals.  For the same rings,
and for right-angled Coxeter groups,
J.~Harlander and H.~Meinert [\ham] have shown that $\vcd_R\Gamma$ is
determined by the local structure of Davis' complex and that Davis'
construction can be generalised to graph products of finite groups.

Our contribution splits into three parts.  Firstly, Davis' complex may
be defined for infinitely generated Coxeter groups (and infinite graph
products of finite groups).  We determine which such groups $\Gamma$ have
finite  virtual cohomological dimension over the integers,  and give
partial information concerning $\vcd_{\Bbb Z}\Gamma$.  We discuss a
form of Poincar\'e duality for simplicial complexes that are like
manifolds from the point of view of $R$-homology, and give conditions
for a (finite-index subgroup of a) Coxeter group to be a Poincar\'e
duality group over $R$.  We give three classes of examples:
we recover Bestvina's examples (and give more information
about their cohomology); we exhibit a group
whose virtual cohomological dimension over the integers is finite but
strictly greater than its virtual cohomological dimension over any
field; we exhibit a torsion-free rational Poincar\'e duality group
which is not an integral Poincar\'e duality group.

Secondly, we discuss presentations for torsion-free subgroups of low
index in right-ang\-led Coxeter groups.  In some cases (depending
on the local structure of Davis' complex) we determine the minimum
number of generators for any torsion-free normal subgroup of minimal
index.  Using the computer package GAP [\gap] we find good
presentations for one of Bestvina's examples, where \lq good' means
having as few generators and relations as possible.

Finally, we give an 8-generator 12-relator presentation of a group
$\Delta$ and a construction of $\Delta$ as a tower of amalgamated free
products, which allows us to describe a good CW-structure for an
Eilenberg-Mac~Lane space $K(\Delta,1)$ and explicitly show that
$\Delta$ has cohomological dimension three over the integers and
cohomological dimension two over the rationals. (In fact $\Delta$ is
isomorphic to a finite-index subgroup of a Coxeter group, but our proofs
do not rely on this.)  The starting point of the work contained in this
paper was the desire to see an explicit Eilenberg-Mac~Lane space for an
example like $\Delta$.

\beginsection 2. Introduction.

A {\it Coxeter system\/} $(\Gamma,V)$ is a group $\Gamma$ and a set of
generators $V$ for $\Gamma$ such that $\Gamma$ has a
presentation of the form
$$\Gamma = \langle V\mid (vw)^{m(v,w)}=1\,\, (v,w\in V) \rangle,$$
where $m(v,v)=1$, and if $v\neq w$ then $m(v,w)=m(w,v)$ is either
an integer greater than or equal to 2, or is infinity (in which
case this relation has no significance and may be omitted).
Note that we do not require that $V$ should be finite.  The
group $\Gamma$ is called a {\it Coxeter group\/}, and in the special
case when each $m(v,w)$ is either 1, 2 or $\infty$, $\Gamma$ is
called a {\it right-angled Coxeter group\/}.

\medskip
\remark Let $(\Gamma,V)$ be a Coxeter system, and let
$m:V\times V\rightarrow \BBb N\cup \{\infty\}$ be the function
occurring in the Coxeter presentation for $\Gamma$.  If $W$
is any subset of $V$ and $\Delta$ the subgroup of $\Gamma$
generated by $W$, then it may be shown that $(\Delta,W)$ is a
Coxeter system, with $m_W$ being the restriction of $m_V$ to
$W\times W$ [\bro].  The function $m$ is determined by
$(\Gamma,V)$ because $m(v,w)$ is the order of $vw$ (which is
half the order of the subgroup
of $\Gamma$ generated by $v$ and $w$).

\medskip
\definition A {\it graph\/} is a 1-dimensional simplicial complex (i.e.\
our graphs contain no loops or multiple edges).  A {\it labelled
graph\/} is a graph with a function from its edge set to a set
of \lq labels'.  A {\it morphism of graphs\/} is a simplicial map which
does not collapse any edges.  A {\it morphism of labelled graphs\/} is
a graph morphism such that the image of each edge is an edge
having the same label.  A {\it colouring\/} of a graph $X$ is a function from
its vertex set to a set of \lq colours' such that the two ends
of any edge have different images.  Colourings of a graph $X$ with
colour set $C$ are in 1-1-correspondence with graph morphisms
from $X$ to the complete graph with vertex set $C$.

\headline{\hbox{\rm Dicks and Leary}
\hss\hbox{\rm On subgroups of Coxeter groups}}
\medskip
\definition For a Coxeter system $(\Gamma,V)$, the simplicial
complex $K(\Gamma,V)$ is defined to have as $n$-simplices the
$(n+1)$-element subsets of $V$ that generate finite subgroups of
$\Gamma$.  Note that our $K(\Gamma,V)$ is Davis' $K_0(\Gamma,V)$
in [\dav].  The graph $K^1(\Gamma,V)$ is by definition the
1-skeleton of this complex.  The graph $K^1(\Gamma,V)$ has a
labelling with labels the integers greater than or equal to 2,
which takes the edge $\{v,w\}$ to $m(v,w)$.
This labelled graph is different from,
but carries the same information as, the Coxeter diagram. The labelled
graph  $K^1(\Gamma,V)$ determines the Coxeter system $(\Gamma,V)$ up to
isomorphism, and any graph labelled by the integers greater than or equal
to 2 may arise in this way.  A morphism of labelled  graphs from
$K^1(\Gamma,V)$ to $K^1(\Delta,W)$ gives rise to a  group homomorphism
from $\Gamma$ to $\Delta$.

Call a subgroup of $\Gamma$ {\it special\/} if it is generated by a
(possibly empty) subset of $V$.  Thus the simplices of
$K(\Gamma,V)$ are in bijective correspondence with the
non-trivial finite special subgroups of $\Gamma$.  Let
$D(\Gamma,V)$ be the simplicial complex associated to the poset
of (left) cosets of finite special subgroups of $\Gamma$.  By
construction $\Gamma$ acts on $D$, and the stabiliser of each
simplex is conjugate to a finite special subgroup of $\Gamma$.
In [\dav], Davis showed that $D(\Gamma,V)$ is contractible if
$V$ is finite, and the general case follows easily (for example
because any cycle (resp. based loop) in $D(\Gamma,V)$ is
contained in a subcomplex isomorphic to $D(\langle V'\rangle,V')$
for some finite subset $V'$ of $V$, so {\it a fortiori\/} bounds
(resp. bounds a disc) in $D(\Gamma,V)$).  Note that
$K(\Gamma,V)$ is finite-dimensional if and only if $D(\Gamma,V)$
is, and in this case the dimension of $D(\Gamma,V)$ is one more
than the dimension of $K(\Gamma,V)$.

A {\it graph product\/} $\Gamma$ of finite groups
in the sense of E.~R.~Green [\gre]
is the quotient of the free product of a family $\{G_v\mid v\in V\}$
of finite groups by the normal subgroup generated by the sets
$\{[g,h]\mid g\in G_v,\,h\in G_w\}$ for some pairs $v\neq w$
of elements of $V$.  A graph product of groups of
order two is a right-angled Coxeter group.  If a {\it special subgroup\/} of a
graph product is defined to be a subgroup generated by some subset of
the given family of finite groups, then the above definitions of
$K(\Gamma,V)$ and $D(\Gamma,V)$ go through unchanged.  In [\ham] it is
proved that for a graph product of finite groups, $D(\Gamma,V)$ is
contractible.  (As in [\dav] only the case when $V$ is
finite is considered, but the general case follows easily.)
The group algebra for a graph product $\Zz\Gamma$ is
isomorphic to the quotient of the free coproduct of the $\Zz G_v$ by
relations that ensure that the pairs $\Zz G_v$ and $\Zz G_w$ generate
 their tensor product whenever
$G_v$ and $G_w$ commute.
Theorem~4.1 of [\dil] is a result for algebras formed in this way
which in the case of the group algebra of a graph product is
equivalent to the acyclicity of $D(\Gamma,V)$.

\beginsection 3. Virtual cohomology of Coxeter groups.

Henceforth we shall make use of the abbreviations $\vcd$ and $\cd$
to denote the phrases \lq virtual cohomological dimension\rq and
\lq cohomological dimension\rq, respectively, and when no ring is
specified, these dimensions are understood to be over the ring of
integers.

\proclaim Theorem 1.  The Coxeter group $\Gamma$ has finite
$\vcd$ if and only if there is a labelled graph
morphism from  $K^1(\Gamma,V)$ to some finite labelled graph.

\pf The complex $K= K(\Gamma,V)$ has simplices of arbitrarily large
dimension if and only if $V$ contains arbitrarily large finite
subsets generating finite subgroups of $\Gamma$.  In this case
$\Gamma$ cannot have a torsion-free subgroup of finite index,
and there can be no graph morphism from $K^1$ to any finite graph.
Thus we may assume that $K$ and hence also $D$ are finite-dimensional.
Any torsion-free subgroup of $\Gamma$ acts freely on $D$, and so it
remains to show that if $D$ is finite-dimensional then there is
a labelled graph morphism from $K^1$ to a finite graph
if and only if $\Gamma$ has a finite-index torsion-free subgroup.

As remarked above, a morphism from  $K^1(\Gamma,V)$ to $K^1(\Delta,W)$ gives
rise to a group homomorphism from $\Gamma$ to $\Delta$ in an obvious
way.  Moreover, if $v,v'$ have product
of order $m(v,v')$, then so do their images in $W$, because
the edge $\{v,v'\}$ and its image in $K^1(\Delta,W)$ are both
labelled by $m(v,v')$.  Now if $V'$ is a finite subset of $V$
generating a finite subgroup of $\Gamma$, and $W'$ is its image
in $W$, then it follows that $\langle V'\rangle$ and $\langle
W'\rangle$ have identical Coxeter presentations, so are
isomorphic.  Thus a morphism from $K^1(\Gamma,V)$ to
$K^1(\Delta,W)$ gives rise to a homomorphism from $\Gamma$ to
$\Delta$ which is injective on every finite special subgroup of
$\Gamma$.
Now suppose that there is a morphism from $K^1(\Gamma,V)$ to
$K^1(\Delta,W)$ for some finite $W$.  The finitely generated
Coxeter group $\Delta$ has a finite-index torsion-free subgroup
$\Delta_1$, so let $\Gamma_1$ be the inverse image of this
subgroup in $\Gamma$.  Since $\Gamma_1$ intersects any conjugate
of any finite special subgroup trivially, it follows that
$\Gamma_1$ acts freely on $D(\Gamma,V)$ and is torsion-free.

Conversely, if $\Gamma$ has a finite-index torsion-free subgroup
$\Gamma_1$, which we may assume to be normal, let $Q$ be the
quotient $\Gamma/\Gamma_1$, and build a labelled graph $X$ with vertices
the elements of $Q$ of order two and all possible edges between
them.  Label the edge $\{q,q'\}$ by the order of
$qq'$.  Now the homomorphism from $\Gamma$ onto
$Q$ induces a simplicial map from $K^1(\Gamma,V)$ to $X$ which
is a labelled graph morphism because if $vv'$ has finite order then its
image in $Q$ has the same order.   \qed

\remark 1) If we are interested only in right-angled Coxeter
groups then all the edges of $K^1$ have the same label, 2,
and we may replace the condition that there is a morphism from
$K^1$ to a finite labelled graph by the equivalent condition that $K^1$
admits a finite colouring.
The above proof can be simplified slightly in this case, because
the right-angled Coxeter group corresponding to a finite
complete graph is a finite direct product of cyclic groups of order
two.

2) An easy modification of the proof of Theorem~1 shows that a graph
product $\Gamma$ of finite groups has finite $\vcd$ if and only if
there are only finitely many isomorphism types among the vertex groups
$G_v$, and the graph $K^1(\Gamma,V)$ admits a finite colouring.

3) Let $(\Gamma,V)$ be the Coxeter system corresponding to the
complete graph on an infinite set, where each edge is labelled
$n$ for some fixed $n\geq 3$.  Then $K(\Gamma,V)$ is
one-dimensional (since the Coxeter group on three generators such that
the product of any two has order $n$ is infinite), $\Gamma$ has an
action on a 2-dimensional contractible complex with stabilisers
of orders 1, 2, and $2n$, but by the above theorem $\Gamma$ does
not have finite $\vcd$.  Similarly, if we take a triangle-free
graph which cannot be finitely coloured, then the corresponding
right-angled Coxeter group acts on a contractible 2-dimensional
complex with stabilisers of orders 1, 2, and 4, but does not
have finite $\vcd$.  In contrast, any group acting on a tree
with finite stabilisers of bounded order has finite $\vcd$; see
for example [\did], Theorem I.7.4.

If a Coxeter group $\Gamma$ has finite $\vcd$ then $D(\Gamma,V)$
is finite-dimensional and the dimension of $D$ gives an upper
bound for $\vcd\Gamma$.  Parts a) and c) of the following theorem
determine when this upper bound is attained.  The information
concerning the right $\Gamma$-module structure on various cohomology
groups will be used only during the construction (in example~3 of
the next section) of a torsion-free rational Poincar\'e duality
group that is not a Poincar\'e duality group over the integers.
To avoid cluttering the statement unnecessarily we first
give some definitions that are used in it.
\medskip
\definition For a Coxeter system $(\Gamma,V)$ and an abelian
group~$A$, let $A^\circ$ denote the
$\Gamma$-bimodule with underlying
additive group $A$ and $\Gamma$-actions given by $va=av=-a$ for all
$v\in V$.  This does define compatible actions
of~$\Gamma$ because each of
the relators in the Coxeter presentation for $\Gamma$ has even
length as a word in $V$.  For a $\Gamma$-module~$M$, let $M_a$
denote the underlying abelian group.  For a simplicial complex~$D$,
let $C_*(D)$ denote the simplicial chain complex of $D$,
let $C^+_*(D)$ denote the augmented simplicial chain complex (having
a $-1$-simplex equal to the boundary of every 0-simplex) and let
$\tilde H^*(D;A)$ denote the reduced cohomology of $D$ with
coefficients in $A$, i.e.\ the  homology of the cochain complex
$\Hom(C^+_*(D),A)$.  All our  $\Gamma$-modules (in particular, all our
chain complexes of  $\Gamma$-modules) are left modules unless otherwise
stated.

\def\pra{\par}
\proclaim Theorem 2.  Let $(\Gamma,V)$ be a Coxeter system such that
$\Gamma$ has finite $\vcd$, let $K=K(\Gamma,V)$ have dimension
$n$ (which implies that $\vcd\Gamma\leq n+1$),
let $D=D(\Gamma,V)$, let $\Gamma_1$
be a finite-index torsion-free subgroup of $\Gamma$, and let
$A$ be an abelian group containing no elements of order two.  Then
\pra\hangindent=\parindent
a) For any $\Gamma_1$-module $M$, $H^{n+1}(\Gamma_1;M)$ is a
quotient of a finite direct sum of copies of $\tilde
H^n(K;M_a)$.
\pra\hangindent=\parindent
b) For each $j$, there is an isomorphism of right
$\Gamma$-modules as follows.
$$H^{j+1}\Hom_\Gamma(C_*(D),A^\circ)\cong\tilde H^j(K;A)^\circ$$
\pra\hangindent=\parindent
c) The right $\Gamma$-module $H^{n+1}(\Gamma;A\Gamma)$
(which is isomorphic to $H^{n+1}(\Gamma_1;A\Gamma_1)$ as a right
$\Gamma_1$-module)  admits a surjective homomorphism onto $\tilde
H^n(K;A)^\circ$.   \pra\hangindent=\parindent
d) If multiplication by the order of each finite special subgroup
of $\Gamma$ induces an isomorphism of $A$, then for each $j$ the right
$\Gamma$-modules $H^{j+1}(\Gamma;A^\circ)$ and $\tilde H^j(K;A)^\circ$
are isomorphic.

\medskip
\pf Let $K'$ be the simplicial complex associated to the poset
of non-trivial finite special subgroups of $\Gamma$, so that
$K'$ is the barycentric subdivision of $K$.  Let $D'$ be the
complex associated to the poset of cosets of non-trivial finite
special subgroups of $\Gamma$.  Then $D'$ is a subcomplex of
$D$, and consists of all the simplices of $D$ whose stabiliser
is non-trivial.  We obtain a short exact sequence of
chain complexes of $\Zz\Gamma$-modules
$$0\rightarrow C_*(D')\rightarrow C_*(D)
\rightarrow C_*(D,D')\rightarrow 0,\eqno{(*)}$$
such that for each $n$ the corresponding short exact sequence of
$\Zz\Gamma$-modules is split.

There is a chain complex isomorphism as shown below.
$$C_*(D,D')\cong \Zz\Gamma\otimes_\zz C_{*-1}^+(K')$$
Topologically this is because the quotient semi-simplicial
complex $D/D'$ is isomorphic to a wedge of copies of the
suspension of $K'$, with $\Gamma$ acting by permuting the
copies freely and transitively.  More explicitly, one may
identify $m$-simplices of $D$ with equivalence classes of
$(m+2)$-tuples $(\gamma,V_0,\ldots,V_m)$, where $V_0\containedin
\cdots\containedin V_m$ are subsets of $V$ generating finite
subgroups of $\Gamma$, $\gamma$ is an element of $\Gamma$,
and two such expressions $(\gamma,V_0,\ldots,V_m)$ and
$(\gamma',V_0',\ldots,V_m')$ are equivalent if $V_i=V_i'$ for all $i$ and
the cosets $\gamma\langle V_0\rangle$ and $\gamma' \langle
V_0\rangle$ are equal.  A map from $C_*(D)$ to
$\Zz\Gamma\otimes_\zz C^+_{*-1}(K')$ may be defined by
$$(\gamma,V_0,\ldots,V_m)\mapsto\cases{ 0&if $V_0\neq\empty$,\cr
\gamma\otimes(V_1,\ldots,V_m)&if $V_0=\empty$,}$$
and it may be checked that this is a surjective chain map with
kernel $C_*(D')$.

The claim of part a) now follows easily.  Applying
$\Hom_{\Gamma_1}(\,\cdot\,,M)$ to the sequence (*) and taking the
cohomology long exact sequence for this short exact sequence of
cochain complexes, one obtains the following sequence.
$$H^{n+1}\Hom_{\Gamma_1}(C_*(D,D'),M)\rightarrow
H^{n+1}\Hom_{\Gamma_1}(C_*(D),M)\rightarrow 0$$
Now $H^{n+1}\Hom_{\Gamma_1}(C_*(D),M)=H^{n+1}(\Gamma_1;M)$, and
there is a chain of isomorphisms as below.
$$\eqalign{H^{n+1}\Hom_{\Gamma_1}(C_*(D,D'),M)&\cong
H^n\Hom_{\Gamma_1}(\Zz\Gamma\otimes C_*^+(K'),M)\cr
&\cong \bigoplus_{\Gamma/\Gamma_1} H^n\Hom(C_*^+(K'),M_a)\cr
&\cong \bigoplus_{\Gamma/\Gamma_1} \tilde H^n(K;M_a)}$$

To prove b), note that since $A$ has no elements of order two,
there are no non-trivial $\Gamma$-module homomorphisms
from the permutation module $\Zz\Gamma/\langle V'\rangle$ to
$A^\circ$ for any non-empty subset $V'$~of~$V$.  Hence applying
$\Hom_\Gamma(\,\cdot\,,A^\circ)$ to the sequence ($*$) one obtains
an isomorphism of cochain complexes of right $\Gamma$-modules
$$\Hom_\Gamma(C_*(D),A^\circ)\cong\Hom_\Gamma(C_*(D,D'),A^\circ).$$
Taking homology gives the following chain of isomorphisms.
$$\eqalign{H^{j+1}\Hom_{\zz\Gamma}(C_*(D),A^\circ)
&\cong H^{j+1}\Hom_{\zz\Gamma}(C_*(D,D'),A^\circ)\cr
&\cong H^j\Hom_{\zz\Gamma}(\Zz\Gamma\otimes C_*^+(K'),A^\circ)\cr
&\cong H^j\Hom_\zz(C_*^+(K'),A^\circ)\cr
&\cong \tilde H^j(K;A)^\circ}$$

Now d) follows easily.  Let $R$ be the subring of $\Qq$ generated by
the inverses of the orders of the finite special subgroups of
$\Gamma$.  Now $\Hom_{R\Gamma}(R\otimes C_*(D),A^\circ)$ is isomorphic
to $\Hom_{\zz\Gamma}(C_*(D),A^\circ)$, and $R\otimes C_*(D)$ is a
projective resolution for $R$ over $R\Gamma$, so d) follows from b).

For c), note that there is an equivalence of functors
(defined on $\Gamma$-modules) between
$\Hom_\Gamma(\,\cdot\,,A\Gamma)$
and $\Hom_{\Gamma_1}(\,\cdot\,,A\Gamma_1)$.  In particular,
$H^*(\Gamma_1;A\Gamma_1)$ and $H^*(\Gamma;A\Gamma)$ are both
isomorphic to the homology of the
cochain complex $\Hom_\Gamma(C_*(D),A\Gamma)$.

There is a $\Gamma$-bimodule map $\phi$ from $A\Gamma$ to $A^\circ$
sending $a.w$ to $(-1)^la$, where $w$ is any element of $\Gamma$
representable by a word of length $l$ in the elements of $V$.
Consider the following commutative diagram of cochain complexes, where
the vertical maps are induced by $\phi$:
$$\matrix{\Hom_{\Gamma}(C_*(D),A\Gamma)&\rightarrow&
\Hom_{\Gamma}(C_*(D,D'),A\Gamma)\cr
\downarrow&&\downarrow\cr
\Hom_{\Gamma}(C_*(D),A^\circ)&\rightarrow&
\Hom_{\Gamma}(C_*(D,D'),A^\circ).\cr}$$
The horizontal maps are surjective because $C_i(D,D')$ is a
direct summand of $C_i(D)$ for each $i$, and the lower horizontal map
is an isomorphism as in the proof of b).  The right-hand vertical map
is surjective because $C_*(D,D')$ is $\Zz\Gamma$-free, and hence the
left-hand vertical map is surjective.

Since each of the cochain complexes is trivial in degrees greater than
$n+1$, one obtains a surjection
$$H^{n+1}\Hom_{\Gamma}(C_*(D),A\Gamma)\rightarrow
H^{n+1}\Hom_{\Gamma}(C_*(D),A^\circ),$$
and hence by b) a surjection of right $\Gamma$-modules
from $H^{n+1}(\Gamma_1;A\Gamma_1)$ to $\tilde H^n(K;A)^\circ$.
\qed

\remark Parts b) and d) of Theorem 2 do not generalise easily
to graph products of finite groups having finite virtual cohomological
dimension, and we have no application for these statements except
in the Coxeter group case.  We outline the generalisation of
a)~and a weaker version of~c) below.

The statement and proof of part~a) carry
over verbatim, and there is a generalisation of part~c).
If $\Gamma$ is a graph product of finite groups with~$l$
distinct isomorphism types of vertex group such that the graph
$K^1(\Gamma,V)$ can be $m$-coloured, then the graph product
version of Theorem~1 implies that $\Gamma$ admits a finite
quotient $G = G_1\times\cdots\times G_k$ for some $k\leq lm$, where
each $G_i$ is isomorphic to a vertex group of $\Gamma$ and each
finite special subgroup of $\Gamma$ is mapped injectively to $G$
with image of the form $G_{i(1)}\times\cdots\times G_{i(j)}$ for some
subset $\{i(1),\ldots,i(j)\}$ of $\{1,\ldots,k\}$.  Now for
$1\leq i\leq k$, let $x_i\in \Zz G$ be the sum of all the
elements of $G_i$, and let $Z$ be the $\Zz\Gamma$-module defined
as the quotient of $\Zz G$ by the ideal generated by the $x_i$.
This $Z$ is the appropriate generalisation of $\Zz^\circ$ to the
case of a graph product, because it is a quotient of $\Zz\Gamma$
of finite $\Zz$-rank and contains no non-zero element fixed by
any vertex group.  To see this, note that
$$Z\cong\Zz G_1/(x_1)\otimes\cdots\otimes\Zz G_k/(x_k),$$
where $\Gamma$ acts on the $i$th factor via its quotient $G_i$.
Each factor is $\Zz$-free, and the action of $G_i$ on
$\Zz G_i/(x_i)$ has no fixed points, because for example
$\Cc\otimes(\Zz G_i/(x_i))$ does not contain the trivial $\Cc
G_i$-module.

The arguments used in the proof of Theorem~2
may be adapted to prove a statement like that of part~b) for
the module $Z$, namely  that for any $j$,
$$H^{j+1}\Hom_\Gamma(C_*(D),Z)\cong \tilde H^j(K;Z_a).$$ 
 From this it may be deduced that
if $\Gamma_1$ is a torsion-free finite-index subgroup of
$\Gamma$, then $H^{n+1}(\Gamma_1;\Zz\Gamma_1)$
admits $\tilde H^n(K;Z_a)$ as a quotient.
A similar result could then be deduced for any torsion-free
abelian group $A$.  A similar result could also be proved for
$A$ an $\Ff_p$-vector space, for $p$ a prime not dividing the order of
any of the vertex groups, by using the fact that $\Ff_p G$ is semisimple
to deduce that $\Ff_p\otimes Z$ has no fixed points for the action
of any~$G_i$.

\proclaim Corollary~3.  If $(\Gamma,V)$ is a finite Coxeter system
such that the topological realisation $|K|$ of $K=K(\Gamma,V)$ is the
closure of a subspace which is a connected $n$-manifold, then for any
finite-index torsion-free subgroup $\Gamma_1$ of $\Gamma$,
$$H^{n+1}(\Gamma_1;\Zz\Gamma_1)\cong\tilde H^n(K;\Zz).$$

\proof We shall apply the condition on $|K|$ in the following
equivalent form:  Every simplex of the barycentric subdivision
 $K'$ of $K$  is contained in an
$n$-simplex, and any two $n$-simplices of $K'$ may be joined by a path
consisting of alternate $n$-simplices and $(n-1)$-simplices, each
$(n-1)$-simplex being a face of its two neighbours in the path and of no
other $n$-simplex.  It suffices to show that under this hypothesis,
$H^{n+1}(\Gamma_1;\Zz\Gamma_1)$ is a cyclic group, because by
Theorem~2 it admits $\tilde H^n(K;\Zz)$ as a quotient and has the same
exponent as $\tilde H^n(K;\Zz)$.

Recall the description of the $m$-simplices of $D=D(\Gamma,V)$ as
$(m+2)$-tuples as in the proof of Theorem~2.  The boundary of the
simplex $\sigma=(\gamma,V_0,\ldots,V_m)$ is given by
$$d(\sigma)=\sum_{i=0}^m(-1)^i(\gamma,V_0,\ldots,V_{i-1},
V_{i+1},\ldots,V_m),$$
and the action of $\Gamma$ by
$$\gamma'\sigma=(\gamma'\gamma,V_0,\ldots,V_m).$$
The stabiliser of $\sigma$ is $\gamma\langle V_0\rangle \gamma^{-1}$.
In the case when $m=n+1$, $V_i$ must be a subset of $V$ of
cardinality $i$, and $\sigma$ is therefore in a free $\Gamma$-orbit.
For $\sigma$ an $(n+1)$-simplex, define $f_\sigma\in\Hom_\Gamma
(C_{n+1}(D),\Zz\Gamma)$ by the equations
$$f_\sigma(\sigma')=\cases{\gamma'&if $\sigma'=\gamma'\sigma$ for some
$\gamma'\in \Zz\Gamma$,\cr
0&otherwise.\cr}$$
The $f_\sigma$ form a $\Zz$-basis for
$\Hom_\Gamma(C_{n+1}(D),\Zz\Gamma)$, so it suffices to show that for
each $\sigma$ and $\sigma'$, $f_\sigma\pm f_{\sigma'}$ is a
coboundary.

From now on we shall fix
$\sigma=(\gamma,V_0,\ldots,V_{n+1})$, and show that $f_\sigma\pm
f_{\sigma'}$ is a coboundary for various choices of $\sigma'$.  If
$\sigma'=(\gamma,V_0,\ldots,V_{i-1},V_i',V_{i+1},\ldots,V_{n+1})$ for
some $i>0$, let $\tau$ be the $n$-simplex
$(\gamma,V_0,\ldots,V_{i-1},V_{i+1},\ldots,V_{n+1})$.  There are
exactly two $i$-element subsets of $V_{i+1}$ containing $V_{i-1}$, so
$\sigma$ and $\sigma'$ are the only $(n+1)$-simplices of $D$ having
$\tau$ as a face.  Defining $f_\tau$ in the same way as $f_\sigma$ and
$f_{\sigma'}$ (which we can do because $\tau$ is in a free
$\Gamma$-orbit), we see that the coboundary of $f_\tau$ is
$(-1)^i(f_\sigma + f_{\sigma'})$.

If $\sigma'=(\gamma,W_0,\ldots,W_{n+1})$, take a path in $K'$ between
the simplices $(V_1,\ldots,V_{n+1})$ and $(W_1,\ldots,W_{n+1})$ of the
form guaranteed by the hypothesis, and use this to make a similar path
of $(n+1)$- and $n$-simplices between $\sigma$ and $\sigma'$, and use
induction on the length of this path to reduce to the case considered
above.

It will suffice now to consider the case when
$\sigma'=(\gamma\gamma',V_0,\ldots,V_{n+1})$.  By induction on the
length of $\gamma'$ as a word in $V$, it suffices to consider the case
when $\gamma'=v$.  Using the cases done above and the fact that
$\{v\}$ is a vertex of some $n$-simplex of $K'$, we may assume that
$V_1=\{v\}$.  Now the $n$-simplex $\tau=(\gamma,V_1,\ldots,V_{n+1})$
of $D$ is a face of only $\sigma$ and $\sigma'$.  The simplex $\tau$
has stabiliser in $\Gamma$ the subgroup $\gamma\langle v\rangle
\gamma^{-1}$, so we may define
$$g_\tau(\tau')=\cases{\gamma'\gamma(1+v)\gamma^{-1}&if
$\tau'=\gamma'\tau$,\cr
0&otherwise.\cr}$$
It is easy to check that $d(g_\tau)=f_\sigma+f_{\sigma'}$, using the
fact that $f_{\sigma'}(\gamma'\sigma)=\gamma'\gamma v\gamma^{-1}$.
\qed

In the same vein we have the following.

\proclaim Proposition~4.  If $(\Gamma,V)$ is a finite Coxeter system
such that the topological realisation $|K|$ of $K=K(\Gamma,V)$ is the
closure of a subspace which is a connected $n$-manifold,  and
$\Gamma_1$ is a finite-index torsion-free subgroup of $\Gamma$, then
the topological space $|D(\Gamma,V)|/\Gamma_1$ (which is an
Eilenberg-Mac~Lane space for  $\Gamma_1$) is homeomorphic to a CW-complex
with exactly one $(n+1)$-cell.

\proof We shall give only a sketch.  The complex
$D(\Gamma,V)/\Gamma_1$ consists of copies of the cone on
$K(\Gamma,V)'$ indexed by the cosets of $\Gamma_1$ in $\Gamma$, where
each $n$-simplex not containing a cone point is a face of exactly two
$(n+1)$-simplices.  (The simplex $(\Gamma_1\gamma,V_1,\ldots,V_{n+1})$,
where $V_1=\{v\}$, is a face of $(\Gamma_1\gamma,\emptyset,V_1,\ldots,
V_{n+1})$ and $(\Gamma_1\gamma v,\emptyset,V_1,\ldots,V_{n+1})$.)  By
hypothesis and this observation there exists a tree whose vertices
consist of all the $(n+1)$-simplices of $D(\Gamma,V)/\Gamma_1$ and
whose edges are $n$-simplices of $D(\Gamma,V)/\Gamma_1$ which are
faces of exactly two $(n+1)$-simplices.  The ends of an edge of the
tree are of course the two $(n+1)$-simplices containing it.  The union
of the (topological realisations of the) open simplices of such a tree
is homeomorphic to an open $(n+1)$-cell.  The
required CW-complex has $n$-skeleton the simplices of
$D(\Gamma,V)/\Gamma_1$ not in the tree, with a single $(n+1)$-cell
whose interior consists of the union of the open simplices of the
tree.
\qed

\remark  The condition on $K(\Gamma,V)$ occurring in the statements of
Corollary~3 and Proposition~4 is equivalent to \lq $K(\Gamma,V)$ is a
{\it pseudo-manifold\/}' in the sense of [\mun].
Neither Corollary~3 nor Proposition~4 has a good analogue for
graph products, because both rely on the fact that the $n$-simplices
of $D(\Gamma,V)$ in non-free $\Gamma$-orbits are faces of exactly two
$(n+1)$-simplices.

In Theorem~5 we summarize a version of Poincar\'e duality for
simplicial complexes that look like manifolds from the point of view
of $R$-homology for a commutative ring~$R$.  Our treatment is
an extension of that of J.~R.~Munkres book [\mun], which covers
the case when $R=\Zz$.  We also
generalise the account in~[\mun] by allowing a group to act on
our \lq manifolds'.  The proofs are very similar to those in~[\mun]
however, so we shall only sketch them here.
\medskip
\definition Let $R$ be a  commutative ring.  An {\it
$R$-homology $n$-manifold\/} is a locally finite simplicial complex~$L$
such that the link of every $i$-simplex of~$L$ has the same
$R$-homology as an $(n-i-1)$-sphere, where a sphere of negative
dimension is empty.  From this definition it follows that $L$ is an
$n$-dimensional complex, and that every $(n-1)$-simplex of $L$ is a
face of exactly two $n$-simplices.  Thus (the topological realisation of)
every open $(n-1)$-simplex of $L$ has an open neighbourhood
in $|L|$ homeomorphic to a ball in $\BBb R^n$.  Say that $L$ is
{\it orientable\/} if the $n$-simplices of $L$ may be oriented
consistently across every $(n-1)$-simplex.  Call such a choice of
orientations for the $n$-simplices an {\it orientation\/} for~$L$.
If $L$ is connected and
orientable then a choice of orientation for one of the $n$-simplices
of $L$, together with the consistency condition, determines a unique
orientation for $L$.  In particular, a connected $L$ has either two
or zero orientations, and a simply connected $L$ has two.

For any locally finite simplicial complex $L$, the {\it cohomology with
compact supports of $L$ with coefficients in $R$\/}, written
$H_c^*(L;R)$, is the cohomology of the subcomplex of the $R$-valued
simplicial cochains on $L$ consisting of the functions which vanish on
all but finitely many simplices of $L$.  (This graded $R$-submodule
may be defined for any $L$, but is a subcomplex only when $L$ is
locally finite.)

\proclaim Theorem 5.  Fix a  commutative ring $R$, and let
$L$ be a connected $R$-homology $n$-manifold.  Let $\Gamma$ be a group acting
freely and simplicially on $L$.  If $L$ is
orientable, let $R^\circ$ stand for the right $R\Gamma$-module upon
which an element $\gamma$ of $\Gamma$ acts as multiplication by $-1$
if it exchanges the two orientations of $L$ and as the identity if it
preserves the orientations of $L$.  If $R$ has characteristic two,
let $R^\circ$ be $R$ with the trivial right $\Gamma$-action.  Then
if either $L$ is orientable or $R$ has characteristic two, there
is for each $i$ an isomorphism of right $R\Gamma$-modules
$$H_c^i(L;R)\otimes R^\circ\cong H_{n-i}(L;R).$$

\proof The statements and proofs contained in sections 63--65 of
[\mun] hold for $R$-homology manifolds provided that all (co)chain
complexes and (co)homology are taken with coefficients in $R$.
For each simplex $\sigma$ of~$L$, one defines the dual block
$D(\sigma)$ and its boundary exactly as in section~64 of [\mun].
From the point of view of $R$-homology, the dual block to an
$i$-simplex of $L$ looks like an
$(n-i)$-cell, and its boundary looks like the boundary of an
$(n-i)$-cell.  Thus as in Theorem~64.1 of [\mun],
the homology of the dual block complex $D_*(L;R)$ is isomorphic
to the $R$-homology of $L$.  There is a natural bijection
between the dual blocks of $L$ and the simplices of $L$.  This is
clearly preserved by the action of $\Gamma$.
Each choice of orientation on $L$ gives rise to homomorphisms
$$\psi:D_{n-i}(L;R)\otimes_R C_i(L;R)\rightarrow R$$
which behave well with respect to the boundary maps, and have
the property that for any simplex $\sigma$ with dual block
$D(\sigma)$, and simplex $\sigma'$, $\psi(D(\sigma)\otimes\sigma')=\pm
1$ if $\sigma = \sigma'$, and 0 otherwise.  This allows one to identify
$D_{n-*}(L;R)$ with $C_c^*(L;R)$.  With the diagonal action  of $\Gamma$
on $D_{n-i}(L;R)\otimes C_i(L;R)$, $\psi$ is not  $\Gamma$-equivariant,
but gives rise to a $\Gamma$-equivariant map
$$\psi':D_{n-i}(L;R)\otimes_R C_i(L;R)\rightarrow R^\circ,$$ and hence an
$R\Gamma$-isomorphism between $D_{n-*}(L;R)$ and  $C_c^*(L;R)\otimes
R^\circ$.   \qed

\remark The referee pointed out that the sheaf-theoretic proof of
Poincar\'e duality in G.~E.~Bredon's book ([\bre], 207--211) also
affords a proof of Theorem~5.

\proclaim Corollary 6.  Let $R$ be a  commutative ring and
let $L$ be a contractible $R$-ho\-mo\-logy $n$-manifold.
Let $\Gamma$ be a group
and assume that $\Gamma$ admits a free simplicial action on~$L$ with
finitely many orbits of simplices.  Then $\Gamma$ is a Poincar\'e
duality group of dimension~$n$ over~$R$, with orientation module the
module~$R^\circ$ defined in the statement of Theorem~5.
The same result holds if $L$ is assumed only to be orientable and
$R$-acyclic rather than contractible.

\proof The simplicial $R$-chain complex for $L$ is a finite free
$R\Gamma$-resolution for $R$, and hence $\Gamma$ is FP over $R$.  Since
$L$  has only finitely many $\Gamma$-orbits of simplices, the cochain
complexes (of right $R\Gamma$-modules)
$$\Hom_{R\Gamma}(C_*(L),R\Gamma)\quad\hbox{and}\quad
C_c^*(L;R)$$
are isomorphic.  Hence by Theorem~5, the graded right $R\Gamma$-module
$H^*(\Gamma;R\Gamma)$ is isomorphic to $R^\circ$ concentrated in
degree $n$.  Thus $\Gamma$ satisfies condition d) of Definition~V.3.3
of [\did] and is a Poincar\'e duality group as claimed.
\qed

\proclaim Corollary~7.  Let $(\Gamma,V)$ be a Coxeter system, and let
$R$ be a  commutative ring.  If $K(\Gamma,V)$ is an
$R$-homology $n$-sphere (i.e.\ $K(\Gamma,V)$ is an $R$-homology
$n$-manifold and $H_*(K(\Gamma,V);R)$ is isomorphic to the
$R$-homology of an $n$-sphere), then any finite-index torsion-free
subgroup of $\Gamma$ is a Poincar\'e duality group over $R$, of dimension
$n+1$.

\proof It suffices to show that whenever $K=K(\Gamma,V)$ is an
$R$-homology $n$-sphere, $D=D(\Gamma,V)$ is an $R$-homology
$(n+1)$-manifold, because then Corollary~6 may be applied to the
free action of the finite-index torsion-free subgroup of~$\Gamma$
on~$D$.  We shall show that the link of any simplex in $D$ is
isomorphic to either $K'$ or a suspension
(of the correct dimension) of the link of some simplex
in $K'$.  This implies that under the hypothesis on~$K$,
the link of every simplex of $D$ is an $R$-homology $n$-sphere.

Let $\sigma=(\gamma,V_0,\ldots,V_m)$ be an $m$-simplex of $D$, and
without loss of generality we assume that $\gamma=1$.
Then the link of $\sigma$ is the collection of simplices $\sigma'$
of $D$ having no vertex in common with $\sigma$ but such that the
union of the vertex sets of $\sigma$ and $\sigma'$ is the vertex set
of some simplex of $D$.  Thus the link of the simplex $\sigma$ as
above consists of those simplices $(\gamma',U_0,\ldots,U_l)$ of $D$
such that the finite subsets $U_0,\ldots,U_l,V_0,\ldots,V_m$
of~$V$ are all distinct, generate finite subgroups of $\Gamma$,
and are linearly ordered by inclusion, where $\gamma'$ is an element
of the subgroup $\Gamma_0$ of $\Gamma$ generated by $V_0$.
The link of $\sigma$
decomposes as a join of pieces corresponding to posets
of the three types listed below, where we adopt the convention
that the join of a complex $X$ with an empty complex is isomorphic
to $X$, and spheres of dimension $-1$ are empty.
\par
1) The poset of all subsets $U$ of $V$ such that $\langle U\rangle$
is finite and $U$ properly contains $V_m$.  This is isomorphic to the
poset of faces of the link in $K'$ of any simplex of $K'$ of dimension
$|V_m|-1$ of the form $(V'_1,\ldots,V'_i,\ldots,V'_{m'})$, where
$m'=|V_m|$ and
$V'_{m'}=V_m$.  By the hypothesis on $K$, this is an $R$-homology
sphere of dimension $n-|V_m|$.
\par
2) For each $i$ such that $0\leq i< m$, the poset of all subsets of
$V$ properly containing $V_i$ and properly contained in $V_{i+1}$.
This is isomorphic to the poset of faces of the boundary of
a simplex with vertex set $V_{i+1}-V_i$, so is a triangulation
of a sphere of dimension $|V_{i+1}|-|V_i|-2$.
\par
3) The poset of all cosets in $\Gamma_0$ of proper special subgroups
of $(\Gamma_0,V_0)$.  (Recall that we defined $\Gamma_0=\langle V_0\rangle$.)
This is a triangulation of a sphere of dimension $|V_0|-1$ on which
the group $\Gamma_0$ acts with each $v\in V_0$ acting as a reflection
in a hyperplane (see [\bro], I.5, especially I.5H).

The link of the $m$-simplex $\sigma$ consists of a join of one piece of
type 1), $m$ pieces of type 2), and one piece of type 3).  All of these
are spheres except that the piece of type 1) is only an $R$-homology
sphere.  It follows that the link of $\sigma$ is an $R$-homology sphere,
whose dimension is equal to the sum
$$\eqalign{(n-|V_m|)+1+&(|V_m|-|V_{m-1}|-2)+1+\cdots\cr
+1+&(|V_2|-|V_1|-2)+1+(|V_1|-|V_0|-2)+1+(|V_0|-1)=n-m.}$$
This sum is obtained from the fact that the dimension of the
join of two simplicial complexes is equal to the sum of their
dimensions plus one,
which is correct in all cases, provided that the empty complex is
deemed to have dimension equal to~$-1$.
\qed
\medskip
\remark 1) Theorem~10 may also be used to prove Corollary~7.

2) A $\Zz$-homology sphere is a {\it generalised homology
sphere\/} in the sense of [\dav].  Davis shows that if $(\Gamma,V)$ is
a Coxeter group such that $K(\Gamma,V)$ is a manifold and a $\Zz$-homology
sphere,
then $\Gamma$ acts on an acyclic manifold with finite stabilisers
(see Sections 12~and~17 of~[\dav]).

\beginsection 4. Unusual cohomological behaviour.

Not every simplicial complex may be $K(\Gamma,V)$ for some
Coxeter system $(\Gamma,V)$; for example the 2-skeleton of a
6-simplex cannot occur.  To see this, note that any labelling of
the edges of a 6-simplex by the labelling set
$\{\hbox{red, blue}\}$ that contains no red triangle must contain
a vertex incident with at least three blue edges.  Now recall
that the Coxeter group corresponding to a labelled triangle can
be finite only if one of the edges has label two.  Writing the
integer two in blue, and other integers in red, one sees that
any 7-generator Coxeter group with all 3-generator special
subgroups finite has a 3-generator special subgroup which
commutes with a fourth member of the generating set.

A condition equivalent to a complex $K$ being equal
to $K(\Gamma,V)$ for
some right-angled Coxeter system $(\Gamma,V)$ is that whenever
$K$ contains all possible edges between a finite set of
vertices, this set should be the vertex set of some simplex of
$K$.  Complexes satisfying this condition are called \lq full
simplicial complexes' or \lq flag complexes' [\bes], [\bro].
The barycentric subdivision of any complex satisfies this
condition.  The barycentric subdivision of an $n$-dimensional
complex admits a colouring with $n+1$ colours, where the
barycentre $\hat\sigma$ of an $i$-simplex $\sigma$ is given the
colour $i\in\{0,\ldots,n\}$.  This proves the following (see
11.3 of~[\dav]).

\proclaim Proposition~8.  The barycentric subdivision of any
$n$-dimensional simplicial complex is isomorphic to
$K(\Gamma,V)$ for some right-angled Coxeter system $(\Gamma,V)$
such that $\vcd\Gamma$ is finite.  \qed

We refer the reader to [\hst] for a statement of the universal
coefficient theorem and a calculation of the $\Ext$-groups
arising in the following examples.

\medskip \noindent
{\bf Example 1 (Bestvina).} (A group of finite cohomological dimension
over the integers whose rational cohomological dimension
is strictly less than its integral cohomological dimension.)
Let $X$ be the space obtained by
attaching a disc to a circle by wrapping its edge around the
circle $n$ times, so that $H_1(X)\cong \Zz/(n)$ and $H_2(X)=0$.
Now let $(\Gamma,V)$ be any Coxeter system such that
$K(\Gamma,V)$ is a triangulation of $X$.  The generating set $V$
will be finite since $X$ is compact, so any such $\Gamma$ will
have finite $\vcd$.  Now if $\Gamma_1$ is a finite-index
torsion-free subgroup of $\Gamma$, $\cd\Gamma$ is at most 3, and
for any $\Zz\Gamma_1$-module $M$, $H^3(\Gamma_1;M)$ is a quotient
of a finite sum of copies of $H^2(X;M_a)$, which is in turn
isomorphic to $\Ext(\Zz/(n),M_a)$ by the universal coefficient
theorem.  In particular, $nH^3(\Gamma_1;M)=0$ for any $M$, and
$H^3(\Gamma_1;\Zz\Gamma_1)\cong\Ext(\Zz/(n),\Zz)\cong\Zz/(n)$
by Corollary~3.  Note that the methods
used by Bestvina [\bes] and by Harlander and Meinert [\ham]
seem to show only that $H^3(\Gamma_1;\Zz\Gamma_1)$
contains elements of order $p$ for each prime $p$ dividing $n$,
whereas our argument gives elements of order exactly $n$.

\medskip\noindent
{\bf Example 2.} (A group whose
cohomological dimension over the integers is finite but strictly
greater than its cohomological dimension over any field.)
Let $X$ be a 2-dimensional CW-complex with
$H_1(X)\cong\Qq$ and $H_2(X)=0$, for example $X$ could be an
Eilenberg-Mac\-~Lane space $K(\Qq,1)$ built from a sequence
$C_1,C_2,\ldots$ of cylinders, where the end of the $i$th
cylinder is attached to the start of the $(i+1)$st cylinder by a
map of degree $i$.  Now let $(\Gamma,V)$ be a Coxeter system
such that $K(\Gamma,V)$ is  a 2-dimensional simplicial complex
homotopy equivalent to $X$, and $\Gamma$
has finite $\vcd$. We shall show that $\vcd\Gamma = 3$, and that for any
field $\Ff$, $\vcd_\ff\Gamma = 2$.  Let $\Gamma_1$ be a finite-index
torsion-free subgroup of $\Gamma$.  Then $\cd\Gamma_1$ is at most 3, and
for any $M$, $H^3(\Gamma_1;M)$ is a quotient of a finite sum of copies of
$H^2(X;M_a)\cong \Ext(\Qq,M_a)$.  If $M$ is an $\Ff\Gamma_1$-module for
$\Ff$ a field of non-zero characteristic~$p$, then $\Ext(\Qq,M_a)$ is an
abelian group which is both divisible and annihilated by $p$, so is
trivial.  If $M$ is an $\Ff\Gamma_1$-module for $\Ff$ a field of
characteristic zero, then $M_a$ is a divisible abelian group, so is
$\Zz$-injective, and so once again $\Ext(\Qq,M_a)=0$.  On the other hand,
$H^3(\Gamma_1;\Zz\Gamma_1)$ is non-zero, because it admits
$\Ext(\Qq,\Zz)$ as a quotient, and $\Ext(\Qq,\Zz)$ is a
$\Qq$-vector space of uncountable dimension.

The group $\Gamma_1$ requires infinitely many generators, but a
2-generator example may be constructed from $\Gamma_1$ using an
embedding theorem of Higman Neumann and Neumann ([\rob],
Theorem~6.4.7).  They show that any countable group $G$ may be
embedded in a 2-generator group $\widehat G$ constructed as an
HNN-extension with base group the free product of $G$ and a free group
of rank two, and associated subgroups free of infinite rank.  An easy
Mayer-Vietoris argument shows that for any ring $R$,
$$\cd_R G\leq \cd_R\widehat G \leq\max\{2,\cd_R G\}.$$
Thus $\widehat\Gamma_1$ is a 2-generator group with $\cd\widehat\Gamma_1=3$
but $\cd_\ff\widehat\Gamma_1=2$ for any field $\Ff$.  We do not know whether
there is a finitely presented group with this property, but the
referee showed us the following Proposition (see also [\bie], 9.12).  

\proclaim Proposition~9.  Let $G$ be a group of type $FP$.  Then
there is a prime field\/ $\Ff$ such that $\cd_\ff G = \cd G$.

\proof Recall that if $G$ is of type $FP$, then for any ring $R$,
$\cd_RG$ is equal to the maximum $n$ such that $H^n(G;RG)$ is
non-zero, and that if $\cd G=n$, then for any $R$, $H^n(G;RG)$ is
isomorphic to $H^n(G;\Zz G)\otimes R$ ([\brp], p199--203).
Let $M$ stand for $H^n(G;\Zz G)$, where $n=\cd G$.
Since $\Hom_G(P,\Zz G)$ is a finitely generated right
$\Zz G$-module for any finitely generated projective $P$, it follows
that $M$ is a finitely generated right $\Zz G$-module.

If $M\otimes \Ff_p$ is non-zero for some prime $p$ we may take $\Ff =
\Ff_p$.  If not, then $M$ is divisible and hence,
as an abelian group, $M$ is a direct sum of a number of copies of
$\Qq$ and a divisible torsion group ([\rob], Theorem~4.1.5).  If
$M\otimes\Qq$ is non-zero then we may take $\Ff = \Qq$.  It remains
to show that $M$ cannot be a divisible torsion abelian group.
Suppose that this is the case, and let $m_1,\ldots,m_r$ be a
generating set for $M$ as a right $\Zz G$-module.
If $N$ is the least common multiple of
the additive orders of the elements $m_1,\ldots,m_r$, multiplication
by $N$ annihilates $M$, contradicting the divisibility of $M$.
\qed

\medskip\noindent
{\bf Example 3.} (A torsion-free rational Poincar\'e duality group
of dimension four which is not an integral Poincar\'e duality group.)
Fix an odd prime $q$, and let $X$ be a lens space with fundamental
group of order $q$, i.e.\ $X$ is a quotient of the 3-sphere by a
free linear action of the cyclic group of order $q$.  It is easy to
see that $X$ is triangulable.  The homology
groups of $X$ are (in ascending order) $\Zz$, $\Zz/(q)$, $\{0\}$ and
$\Zz$.  From the universal coefficient theorem it is easy to see
that $X$ has the same $R$-homology as the 3-sphere for any
commutative ring $R$ in which $q$ is a unit.  Now let $(\Gamma,V)$
be a right-angled Coxeter system such that $K(\Gamma,V)$ is a
triangulation of $X$.  By Corollary~7, any finite-index torsion-free
subgroup $\Gamma_1$ of $\Gamma$ is a Poincar\'e duality group of
dimension four over any $R$ in which $q$ is a unit.

We claim however, that $\Gamma_1$ is not a Poincar\'e
duality group (or PD-group for short) over
the field $\Ff_q$, which implies that it cannot be a
PD-group over the integers.  Since all finite subgroups
of~$\Gamma$ have order a power of two and $q$ is an odd prime, it
follows from Theorem~V.5.5 of [\did] that $\Gamma_1$ is a PD-group
over $\Ff_q$ if and only if $\Gamma$ is.  We shall assume that
$\Gamma$ is a PD-group over $\Ff_q$ and obtain a contradiction.

Firstly, note that the $\Ff_q$-cohomology groups $H^0,\ldots,H^3$ of
$X$ are all isomorphic to $\Ff_q$.  Now it follows from Theorem~2
part~c) that the right $\Gamma$-module $H^4(\Gamma;\Ff_q\Gamma)$
admits $\Ff_q^\circ$ (as defined just above the statement of
Theorem~2) as a quotient.  Thus $\Gamma$ has cohomological
dimension four over~$\Ff_q$, and if $\Gamma$ is a PD-group
over~$\Ff_q$, its orientation module must be~$\Ff_q^\circ$.
In particular, for any $\Ff_q\Gamma$-module~$M$, there should be
an isomorphism for each~$i$
$$H^i(\Gamma;M)=\Ext^i_{\Ff_q\Gamma}(\Ff_q,M)
\cong\Tor^{\ff_q\Gamma}_{4-i}(\Ff_q^\circ,M).$$
Now consider the case when $M$ is the $\Gamma$-bimodule $\Ff_q^\circ$,
viewed as a left $\Ff_q\Gamma$-module.  There is an $\Ff_q$-algebra
automorphism
$\phi$~of~$\Ff_q\Gamma$ defined by $\phi(v)=-v$ for each $v\in V$,
because the relators in $\Gamma$
have even length as words in $V$.  The $\Gamma$-bimodule
obtained from $\Ff_q^\circ$ by letting $\Ff_q\Gamma$ act via $\phi$ is
the trivial bimodule $\Ff_q$.  Thus for each $i$, $\phi$ induces
an isomorphism
$$\Tor^{\ff_q\Gamma}_i(\Ff_q^\circ,\Ff_q^\circ)\cong
\Tor^{\ff_q\Gamma}_i(\Ff_q,\Ff_q).$$
(Here we are viewing the bimodules as left modules when they appear
as the right-hand argument in $\Tor$, and as right modules when they
appear as the left-hand argument.)
Putting this together with the isomorphism obtained earlier, it
follows that
if $\Gamma$ is a PD-group over~$\Ff_q$, then for each~$i$,
$$H^i(\Gamma;\Ff_q^\circ)=\Ext^i_{\ff_q\Gamma}(\Ff_q,\Ff_q^\circ)
\cong\Tor^{\ff_q\Gamma}_{4-i}(\Ff_q,\Ff_q)=H_{4-i}(\Gamma;\Ff_q).
\eqno{(*)}$$
The cohomology groups $H^0,\ldots,H^4$ of $\Gamma$ with coefficients
in $\Ff_q^\circ$ are calculated in Theorem~2 part~d); as vector spaces
over $\Ff_q$ they have dimensions 0, 0, 1, 1, and 1 respectively.
We claim now that any finitely generated right-angled Coxeter group
is $\Ff_q$-acyclic, i.e.\ its homology with coefficients in the
trivial module $\Ff_q$ is 1-dimensional and concentrated in degree
zero.  This leads to a contradiction because
given the claim, the isomorphism (*)
for $i=2$ or~3 is between a 1-dimensional vector space and a
0-dimensional vector space.  The claim follows from Theorem~4.11 of
[\daj], which is proved using an elegant spectral sequence argument.
It is also possible to provide a direct proof by induction on the
number of generators using the fact (see [\chs] or [\gre]) that a finitely
generated right-angled
Coxeter group which is not a finite 2-group is a free
product with amalgamation of two of its proper special subgroups, and
applying the Mayer-Vietoris sequence.

\medskip
\remark
1) Note that the only properties of the space $X$ used in Example 3 are that
$X$ be a compact manifold which is triangulable (in the  weak sense that
it is homeomorphic to the realisation of some  simplicial complex), and
that  for some rings $R$, $X$ be an $R$-homology sphere,
but that there be a prime field $\Ff_q$ for $q\neq 2$ such that $X$ is
not an $\Ff_q$-homology sphere.  These examples show that being
a GD-group over $R$ (in the sense of [\did], V.3.8) is not equivalent
to being a PD-group over $R$.

2) If $\Gamma$ is a Coxeter group such that $K(\Gamma,V)$ is a
triangulation of 3-dimensional real projective space, then Corollary~7
shows that any torsion-free finite-index subgroup of $\Gamma$ is a
PD-group over any $R$ in which 2 is a unit.  The methods used above
do not show that such a group is not a PD-group over $R$ when $2R\neq
R$ however.  The results of the next section show that this is indeed
the case.

\beginsection 5.  Cohomology of Coxeter groups with free coefficients.

The results of this section were shown to us by the referee,
although we are responsible for the proofs given here.
The main result is Theorem~10, which computes $H^*(\Gamma,R\Gamma)$ for any
Coxeter group $\Gamma$ such that $K(\Gamma,V)$ is a (triangulation of
a) closed compact oriented manifold.  This should be contrasted with
Theorem~2, which applies to any Coxeter group, but gives only partial
information concerning cohomology with free coefficients.
Theorem~10 may be used to prove Corollary~7 and to give an alternative
proof of the existence of Example~3 of the previous section.

Let $K$ be a locally finite simplicial complex, so that the cohomology
of $K$ with compact supports, $H_c^*(K)$ is defined, and suppose that
we are given a sequence
$$K_1\supseteq K_2\supseteq K_3\supseteq\cdots$$
of subcomplexes such that each $K_i$ is cofinite (i.e., each $K_i$
contains all but finitely many simplices of $K$), and the intersection
of the $K_i$'s is trivial.  For $j>i$, the inclusion of the pair
$(K,K_j)$ in $(K,K_i)$ gives a map from $H^*(K,K_i)$ to $H^*(K,K_j)$,
and the direct limit is isomorphic to $H_c^*(K)$:
 $$H_c^*(K)\cong \lim_{\rightarrow}(H^*(K,K_1)\rightarrow H^*(K,K_2)
\rightarrow\cdots).\eqno{(**)}$$
(To see this, it suffices to check that a similar isomorphism holds
at the cochain level.)

\proclaim Theorem~10.  Let $R$ be a ring, let
$(\Gamma,V)$ be a Coxeter system such that
$K(\Gamma,V)$ is a triangulation of a closed compact connected
$n$-manifold $X$, and suppose that either $X$ is orientable or\/ $2R=0$.
Then as right $R\Gamma$-module, the cohomology of\/ $\Gamma$ with
free coefficients is:
$$H^i(\Gamma,R\Gamma)=\cases{0&for $i= 0$ or 1,\cr
H^{i-1}(X;R)\otimes_R R\Gamma &for $2\leq i\leq n$,\cr
R^\circ &for $i=n+1$.}$$
(Here $R^\circ$ is the $\Gamma$-module defined above Theorem~2.)

\proof By hypothesis $D=D(\Gamma,V)$ is a locally finite complex, and so
$H^*(\Gamma,R\Gamma)$ is isomorphic to $H_c^*(D;R)$.  Recall that
the complex $D$ may be built up from a union of a collection
of cones on the barycentric subdivision $K'$ of $K(\Gamma,V)$ indexed
by the elements of~$\Gamma$.  For $\gamma\in \Gamma$, let $C(\gamma)$
be the cone with apex the coset $\gamma\{1\}$.  In terms of the
description of $D$ in the proof of Theorem~2, $C(\gamma)$
consists of all simplices of $D$ which may be represented in the
form $(\gamma,V_0,\dots,V_m)$ for some $m$ and subsets $V_0\ldots,V_m$
of $V$.  For any enumeration $1=\gamma_1,\gamma_2,\gamma_3,\ldots$
of the elements of $\Gamma$, define subcomplexes $E_i$ and $D_i$ of $D$ by
$$E_i = \bigcup_{j\leq i} C(\gamma_j),\qquad
D_i = \bigcup_{j>i} C(\gamma_j).$$

Davis' original proof that $D$ is contractible [\dav] uses the
following argument.  For $W$ a subset of $V$, define $K_\sigma(W)$ to
be the subcomplex of $K' = K(\Gamma,V)'$ consisting of the simplices
having a face in common with $K(\langle W\rangle, W)' \subseteq
K(\Gamma,V)$  and all of their faces.  The simplicial interior,
$\sint K_\sigma(W)$, of $K_\sigma(W)$ (i.e., the union of the
topological realisations of the open simplices of $K_\sigma(W)$ which
are not faces of any simplex of $K'- K_\sigma(W)$) deformation
retracts onto $K(\langle W\rangle, W)$ by a linear homotopy.  In
particular, this interior is contractible if $\langle W\rangle$ is
finite, and it may also be shown that in this case $K_\sigma(W)$ is
itself contractible.  In [\dav] an enumeration of the elements of
$\Gamma$ is given such that for each $i$ there exists $W\subseteq V$
with $\langle W\rangle$ a finite subgroup of $\Gamma$, and
$E_i\cap C(\gamma_{i+1})$ is isomorphic to $K_\sigma(W) \subseteq K'$ by
the restriction of the natural isomorphism $C(\gamma_{i+1})\cong CK'$.
By induction it follows that each $E_i$ is contractible, and hence
that $D$ is.

Throughout the remainder of the proof fix an enumeration of $\Gamma$
as in the previous paragraph.  From ($**$) it follows that
$$H^*(\Gamma;R\Gamma)\cong
H_c^*(D;R)\cong \lim_{\rightarrow}(H^*(D,D_i;R)).$$
Let $F_i= D_i\cap E_i$, which could be thought of as the boundary of
$E_i$.  By excision, $H^*(D,D_i;R)$ is isomorphic to $H^*(E_i,F_i;R)$,
and since $E_i$ is contractible, this is in turn isomorphic to the
reduced cohomology group $\tilde H^{*-1}(F_i;R)$.  So far we have used
none of the conditions on $K(\Gamma,V)$ except that $V$ be finite.

The hypothesis that $K(\Gamma,V)$ be a closed $R$-oriented
$n$-manifold is used to compute the limit of the groups
$\tilde H^{*-1}(F_i;R)$.  Roughly speaking, $F_i$ is a connected sum
of $i$ copies of~$K'$.  More precisely, if $W_i\subseteq V$ is such
that $E_i\cap C(\gamma_{i+1})\cong K_\sigma(W_i)$, then $F_{i+1}$ is
obtained from the complexes $F_i- \sint K_\sigma(W_i)$ and
$K'- \sint K_\sigma(W_i)$ by identifying the two copies of
$K_\sigma(W_i)- \sint K_\sigma(W_i)$.  (We defined the
simplicial interior $\sint L$ of a subcomplex $L$~of~$K'$ to be a
topological space, but it is easy to see how to define a subcomplex
$M- \sint L$  of $M$ for any  $L\subseteq M \subseteq K'$.)
Given our hypotheses on $K(\Gamma,V)$, there are equalities
$$H^j(K'- \sint K_\sigma(W_i);R)=\cases{
H^j(K';R)&for $i\neq n$,\cr
0&for $i=n$,\cr}
$$
while $K_\sigma(W_i)-  \sint K_\sigma(W_i)$ is a homology
$(n-1)$-sphere by Poincar\'e-Lefschetz duality for $K_\sigma(W_i)$.
>From the usual argument used to compute the cohomology of a connected
sum and induction it follows that
$$H^j(F_i;R)\cong\cases{\bigoplus_{k=1}^iH^j(K';R)&for $0<j<n$,\cr
R&for $j=n$.  \cr}$$
Moreover, the map from $\tilde H^j(F_i;R)$ to $\tilde H^j(F_{i+1};R)$
given by
$$\tilde H^j(F_i;R)\cong H^{j+1}(D,D_i;R)\rightarrow
H^{j+1}(D,D_{i+1};R)\cong \tilde H^j(F_{i+1};R)$$
is the inclusion of the first $i$ direct summands for $j<n$ and the
identity for $j=n$.

As a right $\Gamma$-module, $H^{n+1}_c(D;\Zz)$ is isomorphic to
$\Zz^\circ$ by Theorem~2c), and the claim for general $R$ follows by
the universal coefficient theorem.  To verify the claimed right
$\Gamma$-module structure for $H^j_c(D;R)$ for $j\leq n$, note that
for each $i$, the images of $D_i$ under translation by each of
$\gamma_1^{-1}=1,\ldots,\gamma_i^{-1}$ are contained in $D_1$, and
check that the $i$ corresponding maps
$$\tilde H^{j-1}(K';R)\cong H^j(D,D_1;R)\rightarrow
H^j(D,D_i;R)\cong \tilde H^{j-1}(K';R)\oplus
\cdots\oplus\tilde H^{j-1}(K';R)$$
are the inclusions of the $i$ distinct direct summands.
\qed

\remark A group $G$ is a PD-group over $R$ if and only if it is of
type $FP$ over $R$ and $H^*(G;RG)$ is isomorphic to $R$ concentrated
in a single degree.  It follows that if $\Gamma$ is such that
$K(\Gamma,V)$ satisfies the hypotheses of Theorem~10, then a
finite-index torsion-free subgroup of $\Gamma$ is a PD-group over $R$
if and only if $K(\Gamma,V)$ is an $R$-homology sphere.  For example,
if $K(\Gamma,V)$ is a triangulation of 3-dimensional real projective
space, then a finite-index torsion-free subgroup of $\Gamma$ is a
PD-group of dimension four over $R$ if and only if $2R= R$.

\beginsection 6.  Generating sets for torsion-free subgroups.

In this section we shall consider only right-angled Coxeter groups, so
$(\Gamma,V)$ shall be a right-angled Coxeter system, and $K(\Gamma,V)$
a full simplicial complex.  The numerical information that we give
will not be very useful if $V$ is infinite.
Recall from Section~3 that a colouring $c : V \rightarrow W$ of
the graph $K^1(\Gamma,V)$ with colour set $W$ gives rise to a
homomorphism from $\Gamma$ to a product of copies of the cyclic group
$C_2$ indexed by
the elements of $W$, such that the kernel $\Gamma_1$ is torsion-free.

\proclaim Proposition~11.  Let $(\Gamma,V)$ be a right-angled Coxeter
system, and let $c:V\rightarrow W$ be a colouring of $K^1(\Gamma,V)$.
Let $\Gamma_1$ be the kernel of the induced map from $\Gamma$ to $C_2^W$.
For $w,w'\in W$, define $K^1(\Gamma,V)(w,w')$ to be the largest
subgraph of $K^1(\Gamma,V)$ all of whose vertices have colour $w$ or
$w'$.  Let $S$ be the following set of elements of $\Gamma_1$.
$$S=\{vv'\mid v,v'\in V,\quad c(v)=c(v')\}$$
\pra\noindent
i) If $c$ is such that any two colours $w,w'$ are adjacent in
$K^1(\Gamma,V)$, then  $S$ generates $\Gamma_1$ as a normal subgroup of
$\Gamma$.   \pra\noindent
ii) If $c$ is such that for each $w,w'\in W$ the graph
$K^1(\Gamma,V)(w,w')$ is connected, then $S$ generates
$\Gamma_1$ as a group.

\proof Let $Q$ be the quotient of $\Gamma$ by the normal subgroup
generated by $S$.  Then $Q$ has $C_2^W$ as a quotient, and
the images of $v$ and $v'$ in $Q$ are equal if $c(v)=c(v')$, so
$Q$ is generated by a set of elements of order two bijective
with $W$.  If the colours $w$, $w'$ are adjacent in $K^1(\Gamma,V)$,
there exist $v$, $v'$ with $c(v)=w$, $c(v')=w'$ which commute
in $\Gamma$.  Thus under the hypothesis in i), the relations of $Q$
include relations saying that all pairs of generators commute, and
so $Q$ is isomorphic to $C_2^W$.  The hypothesis in ii) implies
that in i), so it remains to prove that when each
$K^1(\Gamma,V)(w,w')$  is connected, the subgroup generated by $S$ is
normal.  For this  it suffices to show that if $u$, $v$, $v'$ are
elements of $V$ with  $c(v)=c(v')$, then $uvv'u$ is in the subgroup
generated by $S$.   Let
$$\hbox{$v=v_1$, $u_1$, $v_2$,
$u_2,\ldots,v_{n-1}$, $u_{n-1}$, $v_n=v'$}$$ be the sequence of
vertices on a path in the graph $K^1(\Gamma,V)(c(v),c(u))$ between
$v$ and $v'$.  Thus $c(v_i)=c(v)$, $c(u_i)=c(u)$, and for all
$1\leq i\leq n-1$, $u_i$ commutes with $v_i$ and $v_{i+1}$.
These commutation relations imply that $uvv'u$ is expressible
as the following word in elements of $S$.
$$ (uu_1)(v_1v_2)(u_1u_2)(v_2v_3)\cdots(v_{n-2}v_{n-1})
(u_{n-2}u_{n-1})(v_{n-1}v_n)(u_{n-1}u)= uvv'u\qed$$
\medskip

\remark  The hypothesis in i)  is not very strong.
Given a colouring of a graph in which there exist colours $w$ and $w'$
which are not adjacent, it is possible to identify the colours $w$ and
$w'$ to produce a new colouring of the same graph using fewer colours.

\proclaim Corollary~12.  Let $L$ be an $n$-dimensional simplicial
complex having $N$ simplices in total,
such that every simplex of $L$ is a face of an $n$-simplex and
$|L|-  |L^{n-2}|$ is connected.  If $(\Gamma,V)$ is such
that $K(\Gamma,V)$ is the barycentric subdivision of $L$, then
$\Gamma$ has a torsion-free normal subgroup $\Gamma_1$ of index $2^{n+1}$,
which may be generated by $N-n-1$ elements.  $\Gamma$ has no
torsion-free subgroups of lower index, and any normal subgroup
of this index requires at least this number of generators.

\proof Vertices of $K(\Gamma,V)$ correspond bijectively with simplices
of $L$, and we may colour $K^1(\Gamma,V)$ with the set
$\{0,\ldots,n\}$ by sending a vertex to the dimension of the
corresponding simplex of $L$.  Let $\Gamma_1$ be the kernel of the
induced homomorphism onto $C_2^{n+1}$.  Even without the extra conditions
on $L$, $\Gamma_1$ is a torsion-free subgroup of $\Gamma$ of index
$2^{n+1}$.  Since the vertices of an $n$-simplex of $K(\Gamma,V)$
generate a subgroup of $\Gamma$ isomorphic to $C_2^{n+1}$, $\Gamma$
cannot have a torsion-free subgroup of lower index.  The
abelianisation of $\Gamma$ is isomorphic to $C_2^N$, and so $\Gamma$
cannot be generated by fewer than $N$ elements.  If $\Gamma_2$ is any
normal subgroup of $\Gamma$ of index $2^{n+1}$, then $\Gamma/\Gamma_2$
can be generated by $n+1$ elements, so $\Gamma_2$ cannot be generated
by fewer than $N-n-1$ elements.

Now we claim that
the extra conditions on $L$ are equivalent to the condition that
for any $i,j\in\{0,\ldots,n\}$, the graph $K^1(\Gamma,V)(i,j)$
(as defined in the statement of Proposition~11) is
connected.  Firstly, $|L|- |L^{n-2}|$ is connected
if and only if $K^1(\Gamma,V)(n,n-1)$ is connected.  Now if
$K^1(\Gamma,V)(n,i)$ is connected, every $i$-simplex of $L$ is a
face of some $n$-simplex.  For the converse, note first that in
the special case when $L$ is a single $n$-simplex,
$K^1(\Gamma,V)(i,j)$ is connected for each $i<j\leq n$.  In the case
when $j=n$ this is trivial, and for general $n$ follows by an easy
induction.  For the case of arbitrary $L$, to find a path between
any two vertices $v$, $v'$ of $K^1(\Gamma,V)(i,j)$, first pick
$n$-simplices $u$, $u'$ of $L$ such that $v$ is a face of $u$ and
$v'$ is a face of $u'$.  Then pick a path in $K^1(\Gamma,V)(n,n-1)$
between $u$ and $u'$.  For each $(n-1)$-simplex occurring on this
path, pick one of its $i$-simplices.  This gives a sequence $w_1,
\ldots,w_m$ of $i$-simplices such that $w_l$ and $w_{l+1}$ are
faces of the same $n$-simplex for each $l$, and similarly for
the pairs $v$, $w_1$  and $v'$, $w_n$.  By the special case already
proved, there are paths in $K^1(\Gamma,V)(i,j)$ between each of
these pairs, which concatenate to give a path from $v$ to $v'$.

  Hence by Proposition~11, $\Gamma_1$ can be generated by
the set of pairs $vv'$, where $v$ and $v'$ correspond to
simplices of $L$ of the same dimension.  If we fix for each
dimension $i$ one generator $v_i$, then $vv'=(v_iv)^{-1}(v_iv')$, so
we really  need only the pairs $v_iv$ to generate $\Gamma_1$, and
there are exactly $N-n-1$ of these.
\qed

\remark
1) The conditions imposed on the simplicial complex $L$ are
satisfied if $L$ is a triangulation of a connected $n$-manifold, or
more generally if $|L|$ is the closure of a subset which is a connected
$n$-manifold.  These conditions cannot be omitted.  Let $L$ be the
simplicial bow tie, consisting of two triangles joined at a point,
and let $\Gamma$, $\Gamma_1$ be the Coxeter group and subgroup of index
8 constructed from $L$ as in Corollary~12.  Note that $|L|-
|L^0|$ is not connected.  Using GAP [\gap] it may be
shown that the abelianisation of $\Gamma_1$ is free of rank~11, and so
$\Gamma_1$ requires at least 11 generators, rather than the 10 which
would suffice if Corollary~12 applied.

2) When (as in Corollary~12) the Coxeter group $\Gamma$ has a
torsion-free normal subgroup $\Gamma_1$ of index equal to the order
of the  largest finite subgroup of $\Gamma$, this torsion-free subgroup
will not usually be  unique.

\medskip
There are other ways to construct a torsion-free group having similar
homological properties to a given right-angled Coxeter group.  One
generalisation of the construction given above is as follows.  Given
a right-angled Coxeter group $\Gamma$, and a homomorphism $\psi$ from
$\Gamma$ onto a finite group $Q$ with torsion-free kernel $\Gamma_1$,
let $\Delta$ be any torsion-free group
and $\phi$ any homomorphism from $\Delta$ to $Q$.
The group $P$ defined as the pullback of the following diagram is
torsion-free and has finite index in $\Gamma\times\Delta$.
$$\matrix{P&\mapright{}&\Delta\vphantom{\Delta_f}\cr
\mapdown{}&&\mapdown{\phi}\cr
\Gamma&\mapright{\psi}&Q\vphantom{Q^d}}$$
The group $\Gamma_1$ occurs as such a pullback in the case when
$\Delta$ is the trivial group.  The point about taking different
choices of $\Delta$ is that the resulting group may have a simpler
presentation than $\Gamma_1$.  One such result is the following.

\proclaim Proposition~13.  Let $(\Gamma,V)$ be a right-angled Coxeter
system, and fix an $(n+1)$-colouring $c$ of $K^1(\Gamma,V)$.  Let
$K(\Gamma,V)$ have $M$ edges, and $N$ vertices.  Suppose that the
colouring of $K^1(\Gamma,V)$ has the property that for every $v\in V$,
the star of $v$ contains vertices of all colours.  Then there is a
torsion-free group, $P$, of finite index in
$\Gamma\times {\BBb Z}^{n+1}$ having a presentation with $N$ generators
and $M+N-n-1$ relators, all of length four.  Identifying the
generators of $P$ with the set $V$, the relators are the following
words.
\pra
i) For every edge with ends $v$, $v'$ in $K(\Gamma,V)$, the
commutator $[v,v']$.
\pra
ii) For each colour $w\in W$, for some fixed $v_w$ with $c(v_w)=w$, and
for each $v\neq v_w$ such that $c(v)=w$, the word $v^2v_w^{-2}$.

\proof Let $P_1$ be the group presented by the above generators and
relations.  It suffices to show that $P_1$ is isomorphic to
the pullback $P$ in the diagram below, where
$W$ is the $(n+1)$-element set of colours, $\phi$ is the natural
projection and $\psi$ is the homomorphism induced by the colouring
$c$.
$$\matrix{ P&\mapright{}&{\BBb Z}^W\vphantom{Z_p}\cr
\mapdown{}&&\mapdown{\phi}\cr
\Gamma&\mapright{\psi}&C_2^W}\vphantom{C^{W^d}}$$
Identify the standard basis for $\BBb Z^W$ with the set $W$.  The
elements $(v,c(v))$ of $\Gamma\times \BBb Z^W$ are naturally bijective
with $V$, and satisfy the relations given in the statement.
Thus there is a homomorphism from $P_1$ to $P\leq \Gamma\times \BBb Z^W$
which sends $v$ to $(v,c(v))$.  It remains
to show that this homomorphism is injective and is onto $P$.  The
relations given between the elements of $V$ suffice to show that each
$v^2$ is central in $P_1$, because given $v$, $v'$, there exists $v''$
such that $c(v)=c(v'')$ and there is an edge in $K(\Gamma,V)$ between
$v'$ and $v''$.  By applying the relations given, one obtains
$[v^2,v']=[(v'')^2,v']=1$.  Let $P_2$ be the subgroup of $P_1$
generated by the elements $v^2$.  It is now easy to see that $P_2$ is
central in $P_1$, and is free abelian of rank $n+1$.
Under the map from $P_1$ to $P$, $P_2$ is mapped
isomorphically to the kernel of the map from $P$ to $\Gamma$.  The
quotient $P_1/P_2$ has the same presentation as $\Gamma$.  It follows
that $P_1$ is mapped isomorphically to $P$.
\qed

\remark The condition on the colouring in the statement of
Proposition~13 is weaker than the condition of part ii) of
Proposition~11.  It could be restated as saying that \lq every
component of each $K(\Gamma,V)(w,w')$ contains vertices of
both colours $w$ and $w'$'.
In the case when $K(\Gamma,V)$ is the barycentric subdivision
of an $n$-dimensional simplicial complex $L$ and the colouring taken
is the usual \lq colouring by dimension', the condition
is equivalent to the statement that every simplex of $L$
should be contained in an $n$-simplex.

\beginsection 7.  Presentations of some of Bestvina's examples.

In this section and the next we shall give some explicit presentations
of groups whose cohomological dimensions differ over $\Zz$ and $\Qq$.
To simplify notation slightly we shall adopt the convention that if
$x$ is a generator in a group presentation, $\xinv$ denotes $x^{-1}$.
Much of the section will be based around the eleven vertex
full triangulation of the projective plane given in figure~1, where of
course vertices and edges around the boundary are to be identified in
pairs.  Proposition~14 shows that this triangulation is minimal in some
sense.

\midinsert
%\centerline{\BoxedEPSF{Hexagon scaled 850}}
\psfig{figure=hexagon,height=2in,bbllx=-150pt,bblly=-200pt,bburx=50pt,bbury=0pt}
\endinsert

\proclaim Proposition~14.  There is no full triangulation of the
projective plane having fewer than 11 vertices.

\proof
The following statements are
either trivial or followed by their proof.  A triangulation of the
projective plane with $N$ vertices has $3(N-1)$ edges and $2(N-1)$
faces.  A full triangulation of a closed 2-manifold can have no vertex
of valency 3 or less.  The only 2-manifold having a full $N$-vertex
triangulation with a vertex of valency $N-1$ is the disc.  The only
closed 2-manifold having a full triangulation with a vertex of
valency $N-2$ is the 2-sphere.  There is no triangulation of the
projective plane having 7, 8, or 9 vertices, each of which has valency
4 or 5. (Write an equation for the numbers of vertices of each valency
and obtain a negative number of vertices of valency 4.)  There are
triangulations of the projective plane having 6 vertices, all
of valency 5, but they are not full.

There is no 9 vertex full triangulation of the projective plane:
Assume that there is such a triangulation.  Then by the
above we know that it has a vertex of valency~6.
This vertex and its neighbours form a hexagon containing twelve
edges.  There are no further edges between the vertices of this
hexagon, so all the remaining twelve edges contain at least one of the
remaining two vertices.  Hence at least one of these two vertices is
joined to each of the boundary vertices of the hexagon, giving an
eight vertex triangulation of the 2-sphere before adding the final
vertex.

Any 10 vertex full triangulation of the projective plane has no vertex
of valency seven:  Assume that there is such a vertex.  Then it and
its neighbours form a heptagon containing 14 edges.  There are 13
other edges, each of which must contain at least one of the other two
vertices.  Thus either one of these vertices is joined to all of the
boundary vertices of the heptagon and this gives a 9 vertex
triangulation of the 2-sphere before adding the final vertex, or
the final two vertices are joined to each other and to six each of
the boundary vertices of the heptagon, in which case the complex
contains a tetrahedron (consisting of the two final vertices and any
two of the boundary vertices of the heptagon adjacent to both of the
final vertices).

Any 10-vertex full triangulation of the projective plane has at least
four vertices of valency six (there are no vertices of valency higher
than 6 by the above, and the sum of the 10 valencies is 54).  If no
two of these are adjacent, then they all have the same set of
neighbours, but now all the other vertices have valency 4 and the
total is wrong.  Hence we may assume that a pattern of
edges as in figure~2a occurs in the triangulation.

\midinsert
%\centerline{\BoxedEPSF{Notproj scaled 450}}
\psfig{figure=notproj,height=2in,bbllx=-350pt,bblly=-250pt,bburx=50pt,bbury=50pt}
\endinsert

All of the vertices are already in the picture, and so the vertices
marked $V$ and $W$ can have no other neighbours.  Hence the
triangulation contains the pattern of edges shown in figure~2b.
With two vertices of valency 4, there
must be at least six vertices of valency 6.  Hence by symmetry
it may be assumed that the vertex $X$ is one such.  The only possible
new neighbours for $X$ are the vertices $Y$ and $Z$.  Adding edges
$XY$ and $XZ$, together with the faces implied by fullness, a
triangulated disc whose boundary consists of four edges is obtained.
There is only one way to close up this surface without adding new
vertices or violating fullness, and this gives a triangulation of
the 2-sphere (as may be seen either by calculating its Euler
characteristic, or simply from the fact that after removing one face
it may be embedded in the plane).
\qed

\remark The smallest triangulation of the projective plane which is a
barycentric subdivision has 31 vertices.

\medskip
It is easy to see that the triangulation of figure~1 cannot be
3-coloured, and also that it has no 4-colourings in which every vertex
has a neighbour of each of the three other colours (see Proposition~13).
It does have a 4-colouring in which all
but one of the vertices has neighbours of three colours, and even in
which all but one of the colour-pair subgraphs (as defined in the
statement of Proposition~11) are connected.  The colouring with vertex
classes
$$\{a,c,e\},\quad\{b,d,f\},\quad\{g,h,j\},\quad\{i,k\}$$
is one such.  Let ($\Gamma$,V) be the right-angled Coxeter system with
$K(\Gamma,V)$ as in figure~1.  The above 4-colouring gives rise to an
index 16 torsion-free subgroup of $\Gamma$.  It is easy to see that
as a normal subgroup this group is generated by the eight elements
$$\hbox{$ac$, $ae$, $bd$, $bf$, $gh$, $gj$, $ik$, $gigi=[g,i]$.}$$
Moreover, the techniques of the proof of Proposition~11 can be used
to show that the subgroup generated by these elements is already
normal, and hence that these eight elements generate an index 16
subgroup $\Gamma_2$ of $\Gamma$.

It is still possible to improve upon $\Gamma_2$.  Let
elements $l$, $m$, and $n$ generate a product of three cyclic groups
of order two, and define a homomorphism $\psi$ from $\Gamma$ to this
group by the following equations.
$$\eqalign{
\psi(a)=\psi(c)&=\psi(e)=l\cr
\psi(g)=\psi(h)&=\psi(j)=n}
\qquad\eqalign{
\psi(b)=\psi(d)&=\psi(f)=m\cr
\psi(i)=\psi(k)&=lmn}$$
The homomorphism $\psi$ maps each of the maximal special subgroups of
$\Gamma$ isomorphically to the group generated by $l$, $m$ and $n$,
and hence its kernel $\Gamma_1$ is a torsion-free normal subgroup of
$\Gamma$ of index eight, which contains the index sixteen
subgroup $\Gamma_2$ of the previous paragraph.  The element $bcgi$ is
in $\Gamma_1$, but is not in  $\Gamma_2$.
However, since both $b$ and $c$ commute
with $g$ and $i$,  its square is $(bcgi)^2=gigi$.  It follows that the
eight elements  $$\hbox{$ac$, $ae$, $bd$, $bf$, $gh$, $gj$, $ik$,
$bcgi$}$$ generate the group $\Gamma_1$.  Any normal subgroup of $\Gamma$
of  index eight requires at least eight generators (see Proposition~11),
so there is a sense in which $\Gamma_1$ is best possible.

The Euler characteristic $\chi(\Gamma)$ may be calculated using
I.~M.~Chiswell's formula [\chis], and then $\chi(\Gamma_1)$ is equal to
$|\Gamma:\Gamma_1|\chi(\Gamma)$.  In fact, for a group having a
finite Eilenberg-Mac~Lane space (such as $\Gamma_1$), this Euler
characteristic is just the usual (topological) Euler characteristic of the
Eilenberg-Mac~Lane space.  Indeed, Chiswell's formula can be obtained by
using the Davis complex to make a finite Eilenberg-Mac~Lane space for a
torsion-free subgroup of finite index in a Coxeter group, and then
dividing by the index.  Since the complex $K(\Gamma,V)$ has 11
vertices, 30 edges, and 20 2-simplices, while $\Gamma_1$ has index
eight, the formula gives
$$\chi(\Gamma_1)=8(1- {11\over 2}+ {30\over 4} - {20\over 8})=4.$$

Using the  computer algebra package GAP [\gap], together with some
adjustments  suggested by  V.~Felsch, we were able to find
the presentation for $\Gamma_1$ given below, which
has 8 generators and 12 relations of total length 70 as words in the
generators.  (Recall our convention that $\xinv$ stands for $x^{-1}$.)
$$\eqalign{\Gamma_1 = \langle s,\,t,\,u,\,v,\,&w,\,x,\,y,\,z\mid
\yinv sy\sinv ,\,
vxv\xinv ,\,
\zinv^2wx^2w,\,
\cr&
x^2wu\winv u,\,
\yinv uz\yinv z\uinv ,\,
u\winv uz^2\winv ,\,
\zinv y\zinv \tinv yt,\,
uzsu\sinv \zinv \,
\cr&\qquad
\tinv z\winv t\winv z,\,
vzt\zinv \vinv t,\,
vw\yinv vy\winv ,\,
\vinv w\zinv svzs\winv
\rangle}$$
As words in the eleven generators of the Coxeter group $\Gamma$
the above generators are:
$$
s=ca,\,
t=db,\,
u=fb,\,
v=ki,\,
w=eigb,\,
x=hbie,\,
y=jg,\,
z=cigb.
$$

We know that $\Gamma_1$ needs eight generators, and
also that
an Eilenberg-Mac~Lane space $K(\Gamma_1,1)$ must have at least one 3-cell
(because $H^3(\Gamma_1;\BBb Z\Gamma_1)$ is non-zero by Theorem~2
or Corollary~3).  We also
know that the Euler characteristic of $\Gamma_1$ is 4, and it follows
that the above presentation has the minimum possible numbers of
generators and relations.  Proposition~4 implies that there is a
$K(\Gamma_1,1)$ of dimension three having exactly one 3-cell, but it
does not follow that one may make a $K(\Gamma_1,1)$ by attaching one
3-cell to the 2-complex for a presentation with 8~generators and
12~relations.

In the next section we shall give another presentation for $\Gamma_1$
(although we shall not prove this), also having the minimum numbers
of generators and relations, but with the total length of the
relations slightly longer than here.  The advantage of the other
presentation is that it shows how the group presented (which is in
fact $\Gamma_1$) can be built up using
free products with amalgamation from surface groups, and gives an
independent proof that the cohomological dimension of
$\Gamma_1$ over a ring $R$ depends on whether 2 is a unit in $R$.
We also describe an attaching map for a 3-cell to make an
Eilenberg-Mac~Lane space $K(\Gamma_1,1)$ from the 2-complex for the
new presentation.

It is worth noting that the technique used in Proposition~13 may
be applied to the group $\Gamma$ to give an 11 generator group of
cohomological dimension five over any ring $R$ such that $2R=R$, and
six over other rings, with thirty-seven relators of length four, and
one relator of length eight.  The relators are the thirty commutators
corresponding to the edges in figure~1, together with the following
words.
$$\hbox{$a^2{\inv c}^2$, $a^2{\inv e}^2$, $b^2{\inv d}^2$,
$b^2{\inv f}^2$, $g^2{\inv h}^2$, $g^2{\inv j}^2$,
$i^2{\inv k}^2$, $b^2c^2g^2{\inv i}^2$}$$
As in Proposition~13, one verifies that the subgroup generated by the
squares of the eleven generators is central and free abelian of rank
three, and that the quotient group is isomorphic to $\Gamma$.

\beginsection 8.  A handmade Eilenberg-Mac Lane space.

\proclaim Theorem~15.  The group $\Delta$ given by the presentation
below has cohomological dimension three over rings $R$
such that $2R\neq R$, and cohomological dimension two over rings $R$ such
that $R=2R$.
$$\eqalign{
\Delta=\langle s,\,t,\,u,\,v,\,&w,\,x,\,y,\,z \mid
\sinv v\sinv tu\tinv u\vinv ,\,
w^2\vinv^2 ,\,
\sinv v\sinv v\winv s\vinv sv\winv ,\,
\cr &
x^2\vinv \sinv v\sinv ,\,
\xinv w\vinv s\vinv w\xinv s,\,
u\vinv w\xinv u\xinv \winv v,\,
\winv yt\vinv yw\vinv t,\,
\cr &\qquad
\yinv wu\winv ytu\tinv ,\,
y^2s\vinv s\vinv ,\,
zsu\vinv su\vinv z,\,
\cr &\qquad\qquad
\zinv v\uinv v\uinv x\winv zx\winv ,\,
s\winv yzs\winv yx\winv zx\winv
\rangle}$$
There is an Eilenberg-Mac~Lane space $K(\Delta,1)$, having only one
3-cell and whose 2-skeleton is given by the above presentation.  The
2-sphere forming the boundary of the 3-cell is formed from the
hemispheres depicted in figures 3a and 3b.

\picturesout

\remark In fact this group is isomorphic to the index eight
subgroup $\Gamma_1$ of
the 11 generator Coxeter group $\Gamma$ of the previous section.
The following function $\phi$ from the
generating set for $\Delta$ given above to $\Gamma$ extends to a
homomorphism from $\Delta$ to $\Gamma$, because the image of each
relator is the identity element.  Moreover, it is easy to see that the
image of $\phi$ is equal to $\Gamma_1$.  We shall not prove that
the kernel of $\phi$ is trivial.
$$\eqalign{\phi(s)&=ca\cr \phi(w)&=gifa}\quad
\eqalign{\phi(t)&=hbia\cr \phi(x)&=fija}\quad
\eqalign{\phi(u)&=akia\cr \phi(y)&=fgie}\quad
\eqalign{\phi(v)&=bgia\cr \phi(z)&=dgkc}$$

\picturesout

{\par\noindent{\bf Proof of Theorem~15.} }We shall build the group
$\Delta$ in stages via a tower of free products with amalgamation,
obtained by applying an algorithm due to
Chiswell [\chs].  We shall use the same argument at each stage to justify
this process, but shall give less detail as the steps become more
complicated.  Let $\Delta_0$ be the group given by the following
presentation.   $$\Delta_0=\langle s,\,t,\,u,\,v \mid  \sinv v\sinv
tu\tinv u\vinv  \rangle$$ $\Delta_0$ is the fundamental group of a closed
non-orientable surface of Euler characteristic $-2$, and the 2-complex
corresponding to the above presentation is a CW-complex homeomorphic to
this surface. In particular $\Delta_0$ is torsion-free, and the 2-complex
corresponding to the given presentation is an Eilenberg-Mac~Lane space
$K(\Delta_0,1)$.

Now let $F_0$ be the free group on generators $w$ and $w'$.  Define
automorphisms $a$ and $g$ of $F_0$ by the equations
$$w^a=w\winv',\quad {w'}^a=\winv',\quad w^g=w'\winv,\quad
{w'}^g=ww'\winv.$$
It is easy to check that $a$ and $g$ have order two and commute with
each other, so that they generate a subgroup of $\Aut(F_0)$ isomorphic
to the direct product of two cyclic groups of order two.  Let $G_0$
be the subgroup of $\Aut(F_0)$ generated by $F_0$, $a$ and $g$, which
is isomorphic to the split extension with kernel $F_0$ and quotient
of order four generated by $a$ and $g$.
Now define a homomorphism $\psi_0$ from $\Delta_0$ to $G_0$ by
$$\psi_0(s)=a,\quad \psi_0(t)=gw, \quad\psi_0(u)=agw,\quad
\psi_0(v)=w.$$
This does define a homomorphism from $\Delta_0$ to $G_0$, because the
image of the relator is the identity element, as shown below.
$$\eqalign{\psi_0(\sinv v\sinv tu\tinv u\vinv )&=
awa(gw)(agw)(\winv g)(agw)\winv \cr
&=awagwg\cr
&=w\winv' w'\winv \cr
&=1}$$

\picturesout

Now the images of $s\vinv s v$, $v^2$ and $\vinv s\vinv s$ under
$\psi_0$ are $w'$, $w^2$ and $\winv w'\winv$ respectively.  These
elements generate a normal subgroup of $F_0$ of index two, which is
therefore free of rank three. (Subgroups of index two are always
normal, but what one does is to check that the subgroup is normal, and
then show that it has index two by calculating the order of the
quotient.)  It follows that the elements $s\vinv s
v$, $v^2$ and $\vinv s\vinv s$ freely generate a free subgroup of
$\Delta_0$, and that this subgroup is mapped isomorphically by
$\psi_0$ to the subgroup of $F_0$ generated by $w'$, $w^2$ and $\winv
w'\winv$.  Hence a free product with amalgamation may be made from
$\Delta_0$ and $F_0$ by taking the free product and adding the
relations $w'=s\vinv sv$, $w^2=v^2$, and $\winv w'\winv=\vinv s\vinv
s$.  This gives a group $\Delta_1$.  Using the first of the three new
relations to eliminate the generator $w'$, it follows that $\Delta_1$
has a presentation as below.
$$\Delta_1=\langle s,\,t,\,u,\,v,\,w \mid
\sinv v\sinv tu\tinv u\vinv ,\,
w^2\vinv^2 ,\,
\sinv v\sinv v\winv s\vinv sv\winv
\rangle$$
Moreover, the 2-complex corresponding to this presentation is a
$K(\Delta_1,1)$.

\picturesout

Now take a free group $F_1$ of rank three generated by elements $x$,
$x'$ and $x''$.  Define an automorphism $f$ of $F_1$ by
$$x^f= x,\quad {x'}^f= x\xinv' x,\quad {x''}^f= x\xinv'' x.$$
It may be seen that $f$ has order two.  Let $G_1$ be the split
extension with kernel $F_1$ and quotient the group of order two
generated by $f$, or equivalently the subgroup of $\Aut(F_1)$
generated by $F_1$ and $f$.
Define a homomorphism $\psi_1$ from $\Delta_1$ to $G_1$ as below.
$$\psi_1(s)=x',\quad \psi_1(t)=xf,\quad
\psi_1(u)=x'',\quad \psi_1(v)=xf,\quad \psi_1(w)=x$$
As before, to check that this does define a homomorphism it suffices
to verify that $\psi_1$ sends each relator to the identity in $G_1$.
The images of $s$, $u$, $s\vinv sv$, $\winv v\sinv v\winv$, and $\winv
v\uinv v\winv$ under $\psi_1$ are $x'$, $x''$, $x^2$, $\xinv x'\xinv$,
and $\xinv x''\xinv$ respectively.  These five elements generate a
normal subgroup of $F_1$ of index two, which is therefore a free group
on five generators.  It follows that the subgroup of $\Delta_1$
generated by $s$, $u$, $s\vinv sv$, $\winv v\sinv v\winv$, and $\winv
v\uinv v\winv$ is freely generated by these elements and is mapped
isomorphically to a free subgroup of $F_1$ by $\psi_1$.  Hence we may
form an amalgamated free product of $\Delta_1$ and $F_1$, identifying
the two five generator free subgroups via $\psi_1$.  Call the
resulting group $\Delta_2$.  After eliminating the generators $x'$ and
$x''$ and the relations $x'=s$, $x''=u$, the group
$\Delta_2$ has the presentation given below.  Once again, because the
amalgamating subgroup is free, the 2-complex for this presentation is
a $K(\Delta_2,1)$.

$$\eqalign{\Delta_2=\langle s,\,t,\,u,\,v,\,&w,\,x\mid
\sinv v\sinv tu\tinv u\vinv ,\,
w^2\vinv^2 ,\,
\cr &
\sinv v\sinv v\winv s\vinv sv\winv ,\,
x^2\vinv \sinv v\sinv ,\,
\cr &\qquad
\xinv w\vinv s\vinv w\xinv s,\,
u\vinv w\xinv u\xinv \winv v
\rangle}$$

\picturesout

Take a free group $F_2$ of rank three with generators $y$, $y'$,
$y''$, and define automorphisms $f$ and $i$ of $F_2$ by the following
equations.
$$\eqalign{y^f&=y\cr y^i&=\yinv}\qquad
\eqalign{{y'}^f&=\yinv\yinv'\yinv\cr {y'}^i&=y'y^2}\qquad
\eqalign{{y''}^f&=\yinv\yinv'\yinv''y'y\cr {y''}^i&=\yinv''}$$

\picturesout

It may be checked that $i$ and $f$ have order two and commute, so
they generate a subgroup of $\Aut(F_2)$ isomorphic to a direct product
of two cyclic groups of order two.  Let $G_2$ be the subgroup of
$\Aut(F_2)$ generated by $F_2$,
$f$ and $i$.  Define a homomorphism $\psi_2$ from $\Delta_2$ to
$G_2$ by checking that the function defined as follows on the
generators sends each relator to the identity in $G_2$.
$$\eqalign{\psi_2(s)&=1\cr \psi_2(v)&=yf\cr}\qquad
\eqalign{\psi_2(t)&=y'yf\cr\psi_2(w)&=y\cr}\qquad
\eqalign{\psi_2(u)&=y''\cr\psi_2(x)&=fi\cr}$$
Now $\psi_2$ sends the elements $t\vinv$, $t\uinv\tinv$, $v\sinv
v\sinv$, $w\tinv v\winv$, and $wu\winv$ to $y'$, $y''$, $y^2$, $yy'y$
and $yy''\yinv$ respectively, and these five elements generate a
normal subgroup of $F_2$ of index two, which they therefore freely
generate.  Hence we may make an amalgamated free product $\Delta_3$
from the free product of $\Delta_2$ and $F_2$ by adding the five
relations $t\vinv=y',\ldots,wu\winv=yy''\yinv$.  After eliminating
the generators $y'$ and $y''$, we obtain the following presentation
for $\Delta_3$, such that the corresponding 2-complex is a
$K(\Delta_3,1)$.
$$\eqalign{\Delta_3=
\langle s,\,t,\,u,\,v,\,&w,\,x,\,y,\,z \mid
\sinv v\sinv tu\tinv u\vinv ,\,
w^2\vinv^2 ,\,
\sinv v\sinv v\winv s\vinv sv\winv ,\,
\cr &
x^2\vinv \sinv v\sinv ,\,
\xinv w\vinv s\vinv w\xinv s,\,
u\vinv w\xinv u\xinv \winv v,\,
\winv yt\vinv yw\vinv t,\,
\cr &\qquad
\yinv wu\winv ytu\tinv ,\,
y^2s\vinv s\vinv
\rangle}$$

Now take a 1-relator group $F_3$, with presentation
$$F_3=\langle z,\,z',\,z''\mid \zinv'z\zinv'z''zz''\rangle.$$
$F_3$ is the fundamental group of the closed nonorientable surface of
Euler characteristic~$-1$.  The 2-complex corresponding
to this presentation is a cellular decomposition of this surface, so
in particular is a $K(F_3,1)$.  Define automorphisms $c$ and $g$ of
$F_3$ by the following equations.
$$\eqalign{z^c&=\zinv'z\zinv'\cr z^g&=z'\zinv z'}\qquad
\eqalign{{z'}^c&=\zinv'\cr {z'}^g&=z'}\qquad
\eqalign{{z''}^c&=z''\cr{z''}^g&=\zinv''}$$
To check that $c$ and $g$ as defined above do extend to homomorphisms
from $F_3$ to itself, note that
$$(\zinv'z\zinv' z''zz'')^c=z'\zinv'z\zinv'z'z''\zinv'z\zinv'z''
=zz''\zinv'z\zinv'z''$$
is equal (in the free group) to a conjugate of the relator, and
similarly $(\zinv'z\zinv' z''zz'')^g$ is equal to a conjugate of the
inverse of the relator.  It is easy to check that $c$ and $g$ define
commuting involutions in $\Aut(F_3)$ (and even on the free group
with the same generating set).  Now let $G_3$ be the split extension
with kernel $F_3$ and quotient the subgroup of $\Aut(F_3)$ generated
by $c$ and $g$. (Since $F_3$ has trivial centre, $G_3$ is isomorphic to
a subgroup of $\Aut(F_3)$.)

Define a homomorphism $\psi_3$ from $\Delta_3$ to $G_3$ by taking the
following specification on the generators, and checking that the image of
each relator is equal to the identity in $G_3$.
$$\eqalign{\psi_3(s)&=c\cr\psi_3(w)&=g}\qquad
\eqalign{\psi_3(t)&=1\cr\psi_3(x)&=z''g}\qquad
\eqalign{\psi_3(u)&=gzc\cr\psi_3(y)&=gz'c}\qquad
\eqalign{\psi_3(v)&=g\cr\phantom{\psi_3(v)}&\phantom{=g}}$$
Let $H$ be the subgroup of $\Delta_3$ generated by $\yinv w \sinv$,
$x\winv$, $v\uinv\sinv v\uinv\sinv$ and $w\xinv u\vinv u\vinv$.
The images of these elements under $\psi_3$ are $z'$, $z''$, $z^2$ and
$zz''\zinv$ respectively.  These four elements of $F_3$ generate a
normal subgroup $\psi_3(H)$ of index two, which turns out to be an orientable
surface group (necessarily of Euler characteristic $-2$).  If we write
$\alpha=z^2$ and $\beta=zz''\zinv$, then $\psi_3(H)$ may be presented as
follows.  $$\psi_3(H)=\langle z',\,z'',\,\alpha,\,\beta\mid
\inv\alpha\beta\alpha z''\zinv'\inv\beta\zinv''z'\rangle$$
We claim that $\psi_3$ restricted to $H$ is injective.
For this it suffices to show that the word in the four generators for
$H$ mapping to the relator in $\psi_3(H)$ is equal to the identity in
$\Delta_3$.  Expressed in terms of
the generators for $\Delta_3$, the word is
$$s u \vinv s u\vinv w \xinv u \vinv u \vinv v\uinv \sinv v \uinv
\sinv x \winv s \winv y v \uinv v \uinv x \winv w \xinv \yinv w \sinv,$$
so after cyclically reducing this word it suffices to show that
in $\Delta_3$,
$$u \vinv s u\vinv w \xinv u \vinv \sinv v \uinv
\sinv x \winv s \winv y v \uinv v \uinv  \yinv w = 1.$$
A Lyndon (or van Kampen) diagram whose
boundary is this word, made from eighteen
2-cells each of which has boundary one of the nine relators in
$\Delta_3$, is shown in figure 3a.
This diagram was found, with too much effort, using the normal form
for elements of free products with amalgamation.
It follows that we may form an amalgamated free product $\Delta$
of $\Delta_3$ and $F_3$ amalgamating
$H$ and $\psi_3(H)$.  The generators
$z'$ and $z''$ may be eliminated using the relations $z'=\yinv w\sinv$
and $z''=x\winv$, and the resulting presentation of $\Delta$
is the one given in the statement.

To make an Eilenberg-Mac~Lane space for $\Delta$ it suffices to
attach a 3-cell to the 2-complex for the above presentation
corresponding to the relator in the amalgamating subgroup.  The
boundary of this 3-cell is made from two discs (Lyndon diagrams
expressing the relator in $H$ in terms of the relators in $\Delta_3$
and the relator in $\psi_3(H)$ in terms of the relator in $F_3$) and a
cylinder (whose sides represent the identification of the generators
of $H$ with their images under $\psi_3$).  One of the discs is
figure~3a.  After eliminating $z'$ and $z''$ as above, the rest of
the sphere is as shown in figure~3b.

To verify the claim concerning $H^3(\Delta;M)$ we use the free
resolution for $\Zz$ over $\Zz\Delta$ given by the cellular chain
complex for the universal cover of the Eilenberg-Mac~Lane space
$K(\Delta,1)$ constructed above.
In
the 1-skeleton of the sphere illustrated in figure~3, choose,
for each of the twelve relators,  an oriented path
 between the base points of the two  occurrences of the relator.   These
paths determine elements $w_1,\ldots,w_{12}$ of  $\Delta$.  Now for any
$\Delta$-module $M$, $H^3(\Delta;M)$ is  isomorphic to $M/IM$, where $I$
is the right ideal of $\Zz\Delta$ generated by the twelve elements $1\pm
w_1,\ldots,1\pm w_{12}$, and  the sign in $1\pm w_i$ is positive if the
$i$th relator appears in  figure~3 with the same orientation each time,
and is negative otherwise.  In figure~3a, the two copies of the third
relator meet at their base points and have the same orientation.  It
follows that  $1+w_3=2$, and hence that $H^3(\Delta;M)$ is a quotient of
$M/2M$.  This completes the proof of the statement.  With a little more
work it may be shown that $I$ is equal to the ideal of $\Zz\Delta$
generated by 2 and the augmentation ideal, which  implies that for any
$M$,  $$H^3(\Delta;M)\cong M_\Delta/2M_\Delta.$$
We leave this as an exercise.
\qed

\beginsection 9. Further questions.

\noindent
1) We also used GAP [\gap] to try to find good presentations of various
other finite-index torsion-free subgroups of Coxeter groups.  The
examples that we tried include:
\par\noindent
a) The index sixteen subgroup $\Gamma_2$ of the right-angled Coxeter
group $\Gamma$ described in section~7.  Comparison of the Euler
characteristic (which is twice that of $\Gamma_1$, or eight) with the
known minimum number of generators, together with the fact that any
$K(\Gamma_2,1)$ must have at least one 3-cell indicate that the
minimum number of relators in a presentation of $\Gamma_2$ must be at
least sixteen.  Using GAP we were  able to get down only to an
8-generator 17-relator  presentation, but by hand (and then
checking the result using GAP) we were able to eliminate one of the
relators.  The sum of the lengths of the relators in our presentation is
152.  We found a CW-complex $K(\Gamma_2,1)$ having eight 1-cells,
sixteen 2-cells and one 3-cell.    \par\noindent
b) Other index eight normal subgroups of the Coxeter group $\Gamma$ of
Section~7.  If $\psi'$ is any homomorphism from
$\Gamma$ onto a product of three cyclic groups of order two which
restricts to an isomorphism on each maximal special subgroup of
$\Gamma$, then the kernel of $\psi'$ is a torsion-free index eight
normal subgroup of $\Gamma$.  One way to create such a $\psi'$ is to
take $\psi$ (the homomorphism given earlier, with kernel $\Gamma_1$)
and modify it slightly.  We did not find a $\psi'$ whose kernel had a
smaller presentation than the one given for $\Gamma_1$.
\par\noindent
c) Take figure 1, remove the vertex $h$ and all edges leaving it, and
add a new edge between vertices $i$ and $k$.  Label the three boundary
edges with the label three, and give all other edges the label two.
This gives a presentation of a ten-generator Coxeter system $(\Delta,V)$
such that $K(\Delta,V)$ is a triangulation of the projective plane.
$\Delta$ has torsion-free normal subgroups of index 24, and clearly
has no torsion-free subgroup of lower index, since it contains
elements of order three and subgroups isomorphic to $C_2^3$.  The
subgroup we looked at required nine generators.
\par\noindent
d) In [\bes] Bestvina pointed out that a finite-index torsion-free
subgroup $\Gamma_1$ of a Coxeter group $\Gamma$ such that $K(\Gamma,V)$
is an acyclic 2-complex would have cohomological dimension two over
any ring, but might not have a 2-dimensional
Eilenberg-Mac~Lane space.  (A famous conjecture of Eilenberg and Ganea
asserts that any group of cohomological dimension two has a
2-dimensional Eilenberg-Mac~Lane space [\brp].)  Let $K$ be the
barycentric subdivision of the acyclic 2-complex having  five vertices,
ten edges corresponding to the ten pairs of vertices, and six pentagonal
faces corresponding to a conjugacy class in $A_5$ of elements of order
five.  The simplicial complex $K$ is full.   If $(\Gamma,V)$ is the
corresponding right-angled Coxeter system, then  the easy extension
of Corollary~12 to  polyhedral complexes  shows that
colouring $K(\Gamma,V)$ by dimension gives rise to a torsion-free
index eight normal subgroup $\Gamma_1$ of $\Gamma$, requiring exactly
eighteen generators.  The complex $K(\Gamma,V)$ has 21 vertices, 80
edges and 60 2-simplices, so by Chiswell's formula [\chis], the  Euler
characteristic of $\Gamma_1$ is   $$\chi(\Gamma_1)=8 ( 1- {21 \over
2}+{80\over 4}-{60\over 8})=24.$$  Hence if a presentation for $\Gamma_1$
could be found having twenty-three more relators than generators, the
corresponding 2-complex would be a $K(\Gamma_1,1)$.  An argument similar
to that sketched in the proof of Proposition~4 shows that there is a
3-dimensional $K(\Gamma_1,1)$ having exactly six 3-cells.  Presentations
of $\Gamma_1$ arising from this complex will have 29 more relators than
generators.  Using GAP we found an 18-generator 52-relator presentation
for $\Gamma_1$, but were unable to  reduce the number of relators any
further.  The problem of whether there exists a $K(\Gamma_1,1)$ with less
than six 3-cells remains open.

\medskip\noindent
2) The groups exhibited in Section~4 Example~2
(whose cohomological dimension over the integers is three and whose
cohomological dimension over any field is two) are finitely generated,
but cannot be $FP$ by Proposition~9.  The question of
whether there can be similar examples which are $FP(2)$ or even
finitely presented remains open.

\medskip\noindent
3) The result proved in Proposition~14 does not really prove
that the examples of Sections 7~and~8 are minimal, even
in the sense of being finite-index subgroups of Coxeter groups
with the least possible number of generators.  It may be
true that there can be no full simplicial complex having
ten vertices or fewer whose highest non-zero homology group is
non-free, which  would suggest that  no right-angled Coxeter group
on less than eleven generators  can have different virtual
cohomological dimensions over $\Zz$ and $\Qq$.

\medskip\noindent
4) Is there a simpler example of a group  whose
cohomological dimension over~$\Zz$ is finite and strictly
greater than its cohomological dimension over $\Qq$
than the group $\Gamma_1\cong\Delta$ given in Sections 7~and~8?
%Using the
%Higman-Neumann-Neumann embedding theorem as in the construction of
%Example~2 of Section~4, it is possible to embed
%any $n$-generator group $G$ in a 2-generator group $\widehat G$
%with $\cd_R \widehat G = \cd_R G$ for all rings such that
%$\cd_R G\geq 2$.  An Eilenberg-Mac~Lane space $K(\widehat G,1)$ can be
%made from a $K(G,1)$ by attaching 3 1-cells and $(n+1)$ 2-cells.
%Applying this construction to the 8-generator presentations of the
%group $\Gamma_1\cong\Delta$ of Sections 7~and~8, we found a
%2-generator 12-relator group whose cohomological dimension over any
%ring is equal to that of $\Delta$.
Applying the Higman-Neumann-Neumann embedding theorem to $\Delta$
we were able to construct a 2-generator 12-relator presentation
of a group whose cohomological dimension over any ring is equal to
that of $\Delta$.  An Euler characteristic argument
shows that this group requires at least 12 relators.  The total length
of the 12 relators we found was 1,130, so this group can hardly be
said to be simpler than $\Delta$.

The two distinct 8-generator 12-relator presentations for
$\Delta$ given in Sections 7~and~8 have short relators
(i.e.\ simple attaching maps for the 2-cells) and a simple attaching
map for the 3-cell respectively.  Is it possible to make a
presentation for $\Delta$ in which each 2-cell occurs only
twice in the boundary of the 3-cell and such that the sum of the
lengths of the relators is smaller than 96 (the sum of the lengths
of the relators in the presentation given in Theorem~15)?

\medskip\noindent
5) (P. H. Kropholler) Can there be a group $\Gamma$ with $\cd\Gamma =
4$ and $\cd_{\Bbb Q}\Gamma = 2$? Notice that taking direct products of
copies of Bestvina's examples gives groups with arbitrary finite
differences  between their cohomological dimensions over $\BBb Z$ and over
$\BBb Q$.

\medskip\noindent
6) What we call an $R$-homology manifold is really a {\it simplicial\/}
$R$-homology manifold.  One could give a similar definition and an analogue
of Theorem~5 and Corollary~6 for general (locally compact Hausdorff)
topological spaces.  It may be the case that any torsion-free
Poincar\'e duality group over~$R$ acts freely cocompactly and
properly discontinuously on an
orientable $R$-acyclic $R$-homology manifold.

\goodbreak

\medskip\noindent
{\bf Acknowledgements. } We thank the referee for his comments on an
earlier version of this paper, and especially for showing us
Proposition~9 and Theorem~10.
We thank Volkmar Felsch for informing us
of how to modify   GAP  to suit our purposes.  We
gratefully acknowledge that this work was generously funded by the
DGICYT, through grants PB90-0719 and PB93-0900 for the
first-named author, and a post-doctoral fellowship held at the Centre de
Recerca Matem\`atica for the second-named author.  The second-named
author gratefully acknowledges the support of the Max Planck Institut
f\"ur Mathematik and the Leibniz Fellowship Scheme
during the period when this paper was being revised.

\beginsection 10. References.

\frenchspacing

\par\noindent
[\bes] M. Bestvina, The virtual cohomological dimension of
Coxeter groups, London Math. Soc. Lecture Notes {\bf 181} (1993), 19--23.

\par\noindent 
[\bie] R. Bieri, Homological dimension of discrete groups, Queen Mary 
College Mathematics Notes, London (1976).  

\par\noindent
[\bre] G. E. Bredon, Sheaf theory, McGraw-Hill (1967).

\par\noindent
[\brp] K. S. Brown, Cohomology of groups, Graduate Texts
in Mathematics {\bf 87}, Springer Verlag (1982).

\par\noindent
[\bro] K. S. Brown, Buildings, Springer Verlag (1989).

\par\noindent
[\chs] I. M. Chiswell, Right-angled Coxeter groups, London Math. Soc.
Lecture Notes {\bf 112} (1986), 297--304.

\par\noindent
[\chis] I. M. Chiswell, The Euler characteristic of graph products and
of Coxeter groups, London Math. Soc. Lecture Notes {\bf 173} (1992),
36--46.

\par\noindent
[\dav] M. W. Davis, Groups generated by reflections, Ann. of Math.
{\bf 117} (1983), 293--324.

\par\noindent
[\daj] M. W. Davis and T. Januszkiewicz, Convex polytopes, Coxeter
orbifolds and torus actions, Duke Math. Journal {\bf 62} (1991),
417--451.

\par\noindent
[\did] W. Dicks and M. J. Dunwoody, Groups acting on graphs,
Cambridge Studies in Advanced Mathematics {\bf 17}, Cambridge Univ.
Press (1989).

\par\noindent
[\dil] W. Dicks and I. J. Leary, Exact sequences for mixed
coproduct/tensor-product ring constructions, Publicacions
Matem\`atiques {\bf 38} (1994), 89--126.

\par\noindent
[\gre] E. R. Green, Graph products of groups, Ph. D. Thesis,
University of Leeds, 1990.

\par\noindent
[\ham] J. Harlander and H. Meinert, Higher generation subgroup sets
and the virtual cohomological dimension of graph products of finite
groups, to appear in J. London Math. Soc.

\par\noindent
[\hst] P. J. Hilton and U. Stammbach, A course in homological algebra,
Graduate Texts in Mathematics {\bf 4}, Springer Verlag (1971).

\par\noindent
[\mun] J. R. Munkres, Elements of algebraic topology, Addison-Wesley
(1984).

\par\noindent
[\rob] D. S. Robinson, A course in the theory of groups, Graduate
Texts in Mathematics {\bf 80}, Springer Verlag (1982).

\par\noindent
[\gap] M. Sch\"onert et al., GAP (Groups, Algorithms and Programming)
Version 3 release 2, Lehrstuhl D f\"ur Mathematik, RWTH, Aachen (1993).

\end